\documentclass[12pt,draft]{amsart}
\usepackage[all]{xy}
\usepackage{a4wide}
\usepackage{color}

\begin{document}
\newtheorem{cor}{Corollary}[section]
\newtheorem{theorem}[cor]{Theorem}
\newtheorem{prop}[cor]{Proposition}
\newtheorem{lemma}[cor]{Lemma}
\theoremstyle{definition}
\newtheorem{defi}[cor]{Definition}
\theoremstyle{remark}
\newtheorem{remark}[cor]{Remark}
\newtheorem{example}[cor]{Example}

\newcommand{\cD}{{\mathcal D}}
\newcommand{\cM}{{\mathcal M}}
\newcommand{\cR}{{\mathcal R}}
\newcommand{\cT}{{\mathcal T}}
\newcommand{\C}{{\mathbb C}}
\newcommand{\N}{{\mathbb N}}
\newcommand{\R}{{\mathbb R}}
\newcommand{\Z}{{\mathbb Z}}
\newcommand{\hti}{\tilde{h}}
\newcommand{\Et}{\tilde{E}}
\newcommand{\Kt}{\tilde{K}}
\newcommand{\Mt}{\tilde{M}}
\newcommand{\Mrt}{\tilde{M_r}}
\newcommand{\dr}{{\partial}}
\newcommand{\kappab}{\overline{\kappa}}
\newcommand{\pib}{\overline{\pi}}
\newcommand{\Sigmab}{\overline{\Sigma}}
\newcommand{\gd}{\dot{g}}
\newcommand{\diff}{\mbox{Diff}}
\newcommand{\dev}{\mbox{dev}}
\newcommand{\devb}{\overline{\mbox{dev}}}
\newcommand{\devt}{\tilde{\mbox{dev}}}
\newcommand{\db}{\overline{\partial}}
\newcommand{\Sigmat}{\tilde{\Sigma}}

\newcommand{\cunc}{{\mathcal C}^\infty_c}
\newcommand{\cun}{{\mathcal C}^\infty}
\newcommand{\dd}{d_D}
\newcommand{\dmin}{d_{\mathrm{min}}}
\newcommand{\dmax}{d_{\mathrm{max}}}
\newcommand{\Dom}{\mathrm{Dom}}
\newcommand{\dn}{d_\nabla}
\newcommand{\ded}{\delta_D}
\newcommand{\delmin}{\delta_{\mathrm{min}}}
\newcommand{\delmax}{\delta_{\mathrm{max}}}
\newcommand{\hmin}{H_{\mathrm{min}}}
\newcommand{\maxi}{\mathrm{max}}
\newcommand{\oL}{\overline{L}}
\newcommand{\oP}{{\overline{P}}}
\newcommand{\Ran}{\mathrm{Ran}}
\newcommand{\s}{\mathfrak{s}}
\newcommand{\tgamma}{\tilde{\gamma}}

\newcommand{\II}{I\hspace{-0.1cm}I}
\newcommand{\III}{I\hspace{-0.1cm}I\hspace{-0.1cm}I}
\newcommand{\mnote}[1]{\marginpar{\tiny\tt #1}}

\title{Quasi-Fuchsian manifolds with particles}
\author{Sergiu Moroianu}
\address{Institutul de Matematic\u{a} al Academiei Rom\^{a}ne\\
P.O. Box 1-764\\
RO-014700 Bucharest, Romania}
\email{moroianu@alum.mit.edu}
\author{Jean-Marc Schlenker}
\address{Institut de Math\'ematiques de Toulouse,
Universit\'e Paul Sabatier\\
31062 Toulouse cedex 9, France}
\thanks{J.-M. S. was partially supported by the A.N.R. programs RepSurf, 
2006-09, ANR-06-BLAN-0311, GeomEinstein, 2006-09, and FOG, 2007-10, 
ANR-07-BLAN-0251-01.}
\email{jmschlenker@gmail.com}
\date{version of \today}
\subjclass[2000]{58J20}

\begin{abstract}
We consider 3-dimensional hyperbolic cone-manifolds which are ``convex
co-compact'' in a natural sense, with cone singularities along infinite lines.
Such singularities are sometimes used by physicists as models for massive 
spinless point particles. We prove an infinitesimal rigidity statement 
when the angles around the
singular lines are less than $\pi$: any infinitesimal deformation
changes either these angles, or the conformal structure at infinity
with marked points corresponding to the endpoints of the singular lines.
Moreover, any small variation of the conformal structure at infinity and
of the singular angles can be achieved by a unique small deformation of
the cone-manifold structure. These results hold also when the singularities are 
along a graph, i.e., for ``interacting particles''.
\end{abstract}
\maketitle

\section{Introduction}

\subsubsection*{Quasi-Fuchsian hyperbolic 3-manifolds.}

Let $M$ be the interior of a compact manifold with boundary. 
A complete hyperbolic metric $g$ on $M$ is {\it convex co-compact} if $M$
contains a compact subset $K$ which is {\it convex}: any geodesic segment $c$
in $(M,g)$ with endpoints in $K$ is contained in $K$. Such convex co-compact
metrics (considered up to isotopy) determine a conformal structure on the
boundary at infinity of $M$ (also considered up to isotopy), i.e., an 
element of the Teichm\"uller space of $\dr M$.
According to a celebrated theorem of Ahlfors and Bers (see e.g., 
\cite{ahlfors-bers,ahlfors}), convex co-compact
metrics are uniquely determined by the induced conformal structure at
infinity, and all conformal structures on $\dr M$ can be achieved in this way.

A topologically
simple but already interesting instance is obtained when $M$ is the product
of a closed surface $\Sigma$ by an interval. The space of 
convex co-compact metrics on
$\Sigma \times \R$, which are called ``quasi-Fuchsian'' metrics, is
parametrized by the product of two copies of 
the Teichm\"uller space $\cT_\Sigma$
of $\Sigma$, one corresponding to each boundary component of $M$. In this
manner the geometry of quasi-Fuchsian manifolds has much to say on the
Teichm\"uller theory of $\Sigma$; among many examples we can mention the fact
that the renormalized volume of quasi-Fuchsian metrics provides a K\"ahler
potential for the Weil-Petersson metric on Teichm\"uller space, see 
\cite{takhtajan-teo}.

\subsubsection*{Teichm\"uller theory with marked points.}

The main motivation here is to extend these ideas by replacing the
Teichm\"uller space $\cT_\Sigma$ of $\Sigma$ by its Teichm\"uller space with
$N$ marked points, $\cT_{\Sigma,N}$, and by attaching to each marked point
an angle in $(0,\pi)$. The quasi-Fuchsian metrics on
$\Sigma\times \R$ are then replaced by hyperbolic metrics with conical
singularities along infinite lines going from one connected component of the
boundary to the other; the marked points on each boundary component are
the endpoints of those infinite lines, and the numbers attached to the
marked points are the angles around the singular lines.
We require that the total angle around each singular
curve is less than $\pi$, a restriction which appears naturally at different
stages. In the limit case where those angles tend to $0$ we obtain
geometrically finite hyperbolic manifolds with rank one cusps.

The main result of this paper is the 
first step one has to take when extending the quasi-Fuchsian theory to
encompass those manifolds with conical singularities along infinite lines: 
we prove a local deformation result, namely that the small deformations of 
the ``quasi-Fuchsian cone-manifolds'' described above are parametrized by
the small variations of the angles at the singular lines and of the 
conformal structures at infinity, marked by the endpoints of the singular
lines. The results are actually stated in a more general context of ``convex
co-compact cone-manifolds'', again with ``particles'' -- 
cone singularities along infinite arcs. Note that some results in this
direction, albeit in special cases of manifolds with finite volume, were
obtained by Weiss \cite{weiss-global}.

Our results actually hold for cone-manifolds with singularities along graphs
which have a finite number of vertices,
still under the condition that the cone angle at each singular curve is less
than $\pi$ (as in \cite{weiss-local,weiss-global}). 
Under this condition the singular graph has valence 3. The
vertices can be understood heuristically as ``interactions'' of ``particles''.

\medskip

We now describe in a more detailed way the content of the paper.

\subsubsection*{Hyperbolic cone-manifolds.}

Hyperbolic cone-manifolds were introduced by Thurston (see
\cite{thurston-notes}). They are basically hyperbolic manifolds which are
singular along a stratified subset. In the special case of
3-dimensional cone-manifolds with a singular set which is a disjoint union of
curves, a simple definition can be used (and is given at the beginning
of section 3). In this case, the behavior of the metric in the neighborhood of
a point of the singular locus is entirely determined by a real 
number, the total
angle around the singularity, which is locally constant on the singular
locus. When the singularity is along a graph, the behavior of the metric
close to the vertices is more complicated. However, under the condition
that the cone angles are less than $\pi$, the valence of the singular
graph is $3$, and it remains true that the cone angles determine
completely a local model of the metric.

Hodgson and Kerckhoff \cite{HK} considered compact such hyperbolic
cone-manifolds, for which the singular set is a disjoint union of closed
curves. They showed that, when the total angle around each singular 
curve is less than  $2\pi$, those manifolds are infinitesimally rigid: 
any non-trivial small deformation
induces a deformation of the complex angle around at least one of the connected
curves in the singular locus. Weiss \cite{weiss-local} showed that the same rigidity
result holds when the singular locus is a graph, under the condition that 
the angles at the singular curves are less than $\pi$.

The rigidity result of Hodgson and Kerckhoff was extended by Bromberg
\cite{bromberg1}, who considered complete, non-compact hyperbolic
cone-manifolds, again with singular locus a disjoint union of closed curves,
but also with some non-singular infinite ends similar to the ends of 
convex co-compact
hyperbolic 3-manifolds. In this more general case, any non-trivial 
infinitesimal deformation of the hyperbolic metric induces a non-trivial
deformation either of the conformal structure at infinity, or of the
angle around at least one of the connected curves in the singular locus.

\subsubsection*{Convex co-compact manifolds with particles.}

We consider in this paper complete hyperbolic cone-manifolds, with singularities along
a disjoint union of open curves, or possibly along a graph. 
The difference with the situation considered
by Hodgson and Kerckhoff \cite{HK} or by Bromberg \cite{bromberg1} is that the
curves in the singular locus are not compact, but are instead complete, with
endpoints on the boundary at infinity. A precise definition is given in
section 3, it includes a description of a neighborhood of the endpoints,
ensuring in particular that two singular curves can not be asymptotic 
to each other.

It might be possible to extend the setting considered here to include 
hyperbolic cone-manifolds of finite volume, with cusps. This more general
setting is left for further investigations.

We will use the following definition of convexity, which is stronger than the
condition of having locally convex boundary.

\begin{defi} \label{df:convex}
Let $M$ be a hyperbolic cone-manifold. A subset $C\subset M$ is {\it convex}
if it is non-empty and any geodesic segment in $M$ with
endpoints in $C$ is contained  in $C$.
\end{defi}

For instance, with this definition, points are not convex --- unless $M$ is
topologically trivial. It follows from
the definition that the intersection of two convex subsets of $M$ is either
empty or convex. We show in the appendix that, when the angles at the singular
curves of $M$ are less than $\pi$ and under some weak topological assumptions
on $M$, any closed geodesic in $M$ is
contained in any convex subset. It follows that the intersection of two
convex subsets of $M$ is convex (unless $\pi_1(M)=0$).

\begin{defi} \label{df:coco}
Let $M$ be a complete, non-compact hyperbolic cone-manifold, with singular locus a
graph with a finite number of vertices. 
$M$ is {\it convex co-compact with particles} 
if the angles at each singular line is less than $\pi$ 
and if $M$ contains a compact subset $C$ which is convex. 
\end{defi}

It follows from this definition, and more precisely from Lemma \ref{lm:ends}
in the appendix, that $M$ is then homeomorphic to the interior of a compact
manifold with boundary which we will call $N$ ($N$ is actually homeomorphic
to the compact convex subset $C$ in the definition). The singular set of $M$
corresponds under the homeomorphism with a graph $\Gamma$ embedded in $N$, 
such that vertices of $\Gamma$ adjacent to only one edge are in the boundary
of $N$.

We will in particular
use the term {\it interacting particles} when the 
singular locus has at least one vertex, and {\it non-interacting particles}
when the singular locus is a disjoint union of curves.

A special case is of
interest to us, although it does not play
a central role here (except in the title).

\begin{defi}
A convex co-compact manifold with particles $M$ is called a {\it quasifuchsian 
manifold with particles}
if it is homeomorphic to $\Sigma\times \R$, where $\Sigma$ is a closed surface,
with the singular locus corresponding to lines $\{ x_i\}\times \R$, for 
$x_1,\cdots, x_n\in \Sigma$.
\end{defi}

We hope that the local rigidity result proved here, along with some compactness
statements that will be stated elsewhere, can be used to extend to quasifuchsian
manifolds with particles certain results which are either classical or known for 
convex co-compact (non-singular) hyperbolic manifolds: a Bers-type theorem on
the possible conformal structures at infinity, and statements on what induced
metrics or measured bending laminations can be prescribed on the boundary of 
the convex core. Proving those results in the general setting of convex 
co-compact hyperbolic cone-manifolds appears to be more difficult.

Given a hyperbolic manifold (which is not necessarily complete, so that this
includes the complement of the singular set in the 
cone-manifolds considered here) there is a basic setting,
recalled in section 2, which can be used to understand its 
infinitesimal deformations. It uses a description of those deformations 
as closed 1-forms with values in a vector bundle of ``local'' Killing fields 
defined on the manifold, called $E$ here, 
an idea going back to Weil \cite{weil} and recently
used for cone-manifolds by Hodgson and Kerckhoff \cite{HK}.

Among those deformations, some do not change the underlying geometry of the
manifold; they are the differentials (with respect to a  natural flat
connection on $E$) of sections of $E$, they are called {\it trivial}
deformations. 

\subsubsection*{Main statements.}

The first result of this paper is an infinitesimal rigidity result, stating
that infinitesimal deformations of one of the cone-manifolds considered here
always induces a infinitesimal variation of one of the ``parameters'': the
conformal structure at infinity, or the angle around the singular curves.

\begin{theorem} \label{tm:rigidity}
Let $(M,g)$ be a convex co-compact manifold with particles. Any
non-trivial infinitesimal deformation of the hyperbolic cone-metric $g$ 
induces a non-trivial deformation of the conformal structure with marked points
at infinity or of one of the angles around the singular lines. 
\end{theorem}

The second, related, result is that the small deformations of these
``parameters'' are actually in one-to-one correspondence with the small 
deformations of the cone-manifolds. Let $\cR(M_r)$ be the representation 
variety of $\pi_1(M_r)$ into $PSL_2(\C)$ and $\rho$ the holonomy representation 
of $M_r$. We call $\cR_{cone}(M_r)$ the subset of representations for which
the holonomy of meridians of the singular curves have no translation 
component, that is, the holonomy of the meridians are rotations. Thus 
$\rho\in \cR_{cone}(M_r)$, and, in the neighborhood of $\rho$, the points
of $\cR_{cone}(M_r)$ are precisely the holonomies of cone-manifolds. 

\begin{theorem} \label{tm:existence}
Let $(M,g)$ be a convex co-compact manifold with particles. Let
$c$ be the conformal structure at infinity, and let $\theta_1, \cdots,
\theta_N\in(0,\pi)$ be the angles around the singular lines. 
In the neighborhood of $\rho$, the quotient of $\cR_{cone}(M_r)$ 
by $PSL(2,\C)$ is parameterized by
small deformations of $c,\theta_1, \cdots, \theta_n$. 
\end{theorem}

Note that these results could be somewhat extended, at the cost of more
complicated statements but without any significant change in the proof; it
should be possible to include singularities along closed curves, still under
the hypothesis that the angles around those singularities are less than $\pi$
(or perhaps even $2\pi$ as in \cite{HK}). On the other hand, the condition that
the angle around the ``open'' singular curves is less than $\pi$ appears to be
necessary, at least it occurs at several distinct points in 
the proof given here, and it also comes up naturally
in other properties of those cone-manifolds ``with singular infinity'' that
will not be treated here (in particular the geometry of the boundary of their 
convex cores).

\subsubsection*{More about the motivations.}

It was mentioned above that the main motivation for our work is the search for
a generalization to manifolds with particles of the classical
result, due to Ahlfors and Bers, describing convex co-compact hyperbolic
metrics in terms of the conformal structure on their boundary at infinity. 

It appears conceivable that a proof of such a statement could follow a
``deformation'' approach: proving that, given the angles around the singular
lines, the natural map sending a cone-manifold
to its conformal structure at infinity (marked by the position of the
endpoints of the singular lines) is a homeomorphism. The topology that one
should consider on the space of cone-manifolds would then come from the
associated representations of the fundamental group of the complement
of the singular locus. Three main difficulties would arise:
\begin{itemize}
\item showing that the map is a local homeomorphism --- this is precisely the
  content of Theorem \ref{tm:existence},
\item showing that the map is proper --- which translates as a compactness
  question for convex co-compact manifolds with particles (see \cite{conebend}),
\item showing that some conformal data have a unique inverse image --- a point
  which appears not to be difficult for some particular values of the angles
(of the form $2\pi/k, k\in \N$) for which finite coverings can be used. 
\end{itemize}
So, given the results presented here, a kind of ``double uniformization''
theorem for manifolds with particles would follow from some
compactness results. Since such statements depend on geometric methods which
are completely different from the more analytic tools used here, we have
decided not to include any developments concerning them, and hope to
treat them in a subsequent work.

This line of arguments also leads to applications to Teichm\"uller theory, 
in particular for the Teichm\"uller space of hyperbolic metrics with cone
singularities of prescribed angle on a closed surface. Having a Bers-type
theorem for quasifuchsian manifolds with particles would make it possible
to use in this context renormalized volume arguments as those used in 
\cite{volume} to recover results of \cite{takhtajan-teo}, and to show that
the natural Weil-Petersson metric on those Teichm\"uller spaces is K\"ahler,
and has the renormalized volume as a K\"ahler potential.

\subsubsection*{The geometry of the convex core.}

It is possible to define the convex core of hyperbolic cone-manifold as 
for non-singular convex co-compact hyperbolic 3-manifolds. This appears
natural but proving it properly leads to some technical considerations which
have been moved to the appendix to keep the main part of the paper focused. 
With respect to the
properties of the convex core, the ``convex co-compact'' manifolds with
particles that we consider here appear to share some important
properties of (non-singular) convex co-compact hyperbolic manifolds
(the hypothesis that the cone angles are less than $\pi$ is relevant here). 
This is
beyond the scope of this paper, however we do need some definitions, since
they will be helpful in the geometric constructions of sections 3 and 4.

\begin{defi}
The {\it convex core} $CC(M)$ is the smallest non-empty subset of $M$ which is convex.
\end{defi}

The term ``convex'' should be understood here as in Definition
\ref{df:convex}, and ``smallest'' is for the inclusion. The existence of
$CC(M)$ is clear as soon as $M$ is non-contractible. Indeed, $M$ itself is
convex, while the intersection of two convex subsets of $M$ is convex 
and contains any closed geodesic of $M$ (see Lemma \ref{lm:closed-geod}).

It follows from this definition that $CC(M)$ is a convex set
without extremal points (outside the singular locus of $M$), and
therefore that the intersection of its boundary with the regular set
in $M$ is a ``pleated surface'' as for (non-singular) quasi-Fuchsian
manifolds (see \cite{thurston-notes}). When the angles around the singular
lines are less than $\pi$, a simple but interesting phenomenon occurs: 
the convex core $CC(M)$ contains all the vertices, and 
its boundary is ``orthogonal'' to the singular locus of $M$,
so that its induced metric is a hyperbolic metric with cone singularities
(at the intersection with the singular locus) of angle equal to the
angle of the corresponding curve of the singular locus. Moreover, still
under the hypothesis that the singular angles are less than $\pi$, 
the support of the 
bending lamination of $\dr CC(M)$ does not contain its intersection
with the singular lines. 

These aspects of the geometry of quasi-Fuchsian cone-manifolds, which will
not be developed much here, are important as motivations since they 
appear to indicate that several interesting questions concerning 
quasi-Fuchsian manifolds can also be asked for quasi-Fuchsian
manifolds with particles as defined here, for instance whether any couple of
hyperbolic metrics with cone singularities of prescribed angles can
be uniquely obtained as the induced metric on the boundary of the
convex core, or whether any couple of ``reasonable'' measured 
laminations, on a surface with some marked points, can be uniquely obtained 
as the bending lamination of the boundary of the convex core (for
non-singular quasi-Fuchsian manifolds, see \cite{bonahon-otal,lecuire}).
Other similar questions concerning domains with smooth boundary can
also be considered (see \cite{hmcb} for the non-singular analog).

\subsubsection*{AdS manifolds and 3d gravity.}

G. Mess \cite{mess} discovered that there is a class of anti-de Sitter
manifolds, sometimes called ``globally hyperbolic maximal compact'' (GHMC),
which is in many ways analogous to the quasi-Fuchsian hyperbolic manifolds.
One such analogy is the fact that the space of GHMC AdS manifolds of 
given topology is parametrized by the product of two copies of 
Teichm\"uller space, and the geometry of the convex core presents striking
similarities with the quasifuchsian case. 

Cone singularities along time-like lines are quite natural in the
context of those AdS manifolds, since they are used in the physics
literature to model point particles. It appears (see \cite{minsurf,cone})
that some properties of hyperbolic and AdS manifolds with cone
singularities along open lines (which are time-like in the AdS case)
are quite parallel. 

{\small
\subsection*{Acknowledgments.} We are grateful to Sylvain Gol\'enia for 
pointing out a gap from a previous version of Section \ref{Sir}, and to 
an anonymous referee for many helpful comments and remarks which lead 
to considerable improvements.} 

\section{Local deformations}\label{local}

We recall the link between infinitesimal deformations of hyperbolic metrics
and the first cohomology group of the bundle of infinitesimal Killing fields.

\subsection{The developing map of a hyperbolic metric.}\label{sec2.1}

Let $M_r$ be a connected $3$-manifold, with a hyperbolic metric $g$ 
(i.e., Riemannian metric with constant sectional curvature $-1$); 
this metric does not have to be complete, we are interested in the 
regular set of a hyperbolic cone-manifold. Each
point $x\in M_r$ has a neighborhood which is isometric to an open subset of
hyperbolic $3$-space $H^3$. This isometry can be extended uniquely to a
local isometry 
from the universal cover $(\Mrt,g)$ to $H^3$, called the {\it developing
map} of $(M_r,g)$. We denote it by $\dev_g$, it is well defined up to
composition on the left by a global isometry of $H^3$.

If $(M_r,g)$ is the regular part of a hyperbolic cone-manifold $M$, 
then $\dev_g$ is defined outside the
singular set of $M$. It is usually not injective.

\subsubsection*{Deformations of hyperbolic metrics.}

Let $\gd$ be a infinitesimal deformation of the hyperbolic metric $g$; $\gd$ is
a section of the bundle of symmetric bilinear forms over $M$. We suppose that
$\gd$ is such that the metric remains hyperbolic, i.e., the infinitesimal
variation of the sectional curvature of $g$ induced by $\gd$ vanishes. 

One obvious way to define such ``hyperbolic'' deformations of $g$ is by
considering the Lie derivative of $g$ under the action of a vector field $u$ on
$M$. We call such infinitesimal deformations {\it trivial}. 

\subsubsection*{The vector field associated to a deformation}

Consider the germ at $t=0$ of 
a smooth $1$-parameter family $(g_t)_{0\leq t<\epsilon}$ 
of hyperbolic metrics on $M$ with $g_0=g$ and 
$(\partial_t g_t)_{t=0}=\gd$. Choose a smooth 
$1$-parameter family of developing maps $\dev_{g_t}$ for the metrics
$g_t$. One  
way to do this is as follows: fix a point $x_0$ in $\Mt$, a point $p_0$ 
in $H^3$ and an isometry $I$ between $T_{x_0}\Mt$ and $T_{p_0}H^3$, then 
there exists a unique $\dev_{g_t}$ with the property
\begin{align*}\dev_{g_t}(x_0)=p_0,&&(\dev_{g_t})_*(x_0)=I.\end{align*}
Any other choice must be of the form 
\[\dev_{g_t}'=a_t\dev_{g_t}\]
for some smooth family $(a_t)_{0\leq t<\epsilon}$ of isometries of $H^3$.

For each $x\in\Mt$ the curve 
\[t\mapsto \dev_g^{-1}(\dev_{g_t}(x))\] 
is well-defined for some 
positive time, in particular it defines a vector at $x$. Denote
by $u$ the vector field on $\Mt$ obtained in this way. 

Let $G$ be the group of deck transformations of $\Mt$. Then $u$ is 
\emph{automorphic} with respect to the action of $G$, in the sense that 
for all $\gamma\in G$, the vector field $\gamma_* u-u$ is Killing (we follow
here the terminology used in \cite{HK}).
Indeed, by definition $\gamma_* g_t=g_t$, so there exists an isometry
$a_\gamma(t)$ of $H^3$ such that 
\[\dev_{g_t}\circ\gamma=a_\gamma(t)\circ\dev_{g_t}.\]
By differentiation at $t=0$ this implies
\[u_{\gamma x}=\gamma_*u_x+\dev_g^*\dot{a}_\gamma.\]

\subsection{The bundle $E$ of germs of Killing fields}

Over an arbitrary Riemannian manifold $M$ consider the vector bundle
\[E:=TM\oplus\Lambda^2 T^*M\]
with connection $D$ given by 
\[D_V(u,\alpha)=(\nabla_V u+V\lrcorner\alpha, \nabla_V\alpha-R_{uV})\]
where $R$ is the curvature tensor (we identify vectors and $1$-forms 
using the Riemannian metric). Define a differential operator
$\s:\cun(M,TM)\to\cun(M,E)$ by the formula
\[u\mapsto \s_u:=\left(u, -\frac12(\nabla u)_{\text{anti-sym}}\right).\]
The operator $s$ is called the \emph{canonical lift}, see \cite{HK}.
The following elementary lemma is well-known (see \cite{kostant}):
\begin{lemma}\label{lkf}
On every Riemannian manifold, the canonical lift operator
induces an isomorphism between the space of Killing vector fields
and the space of parallel sections of $E$.
\end{lemma}

We specialize now
to $M$ orientable of dimension $3$, so we identify $\Lambda^2T^*M$ with $TM$ via
the Hodge star and duality. Keeping into account that the sectional 
curvature is $-1$, $R_{uV}$ is mapped under this
identification to $u\times V$. Let $v$ be the vector corresponding to the
$2$-form $\alpha$, then $V\lrcorner\alpha$ is dual to $v\times V$. Hence
under this identification of $E$ with $TM\oplus TM$, the connection $D$ becomes
\[D_V(u,v)=(\nabla_Vu+v\times V, \nabla_Vv-u\times V).\]
To simplify even further, note that $E$ is 
isomorphic to the complexified tangent bundle $T_\C M$ via
\[(u,v)\mapsto u+iv.\]
We extend by linearity the Levi-Civit\`a connection and the vector product
to $T_\C M$. Hence the bundle with connection $(E,D)$ is isomorphic to
$T_\C M$ with the connection (again denoted by $D$) given by
\begin{equation}\label{D}
D_V \phi=\nabla_V\phi +iV\times \phi~.
\end{equation}
Clearly $D$ commutes with complex multiplication.
In this framework, the canonical lift operator is given
by the expression
\begin{equation}\label{canlif}
\s_u=u-i\frac{\mathrm{curl} (u)}{2}.
\end{equation}

Using the fact that $M$ is hyperbolic, a straightforward computation 
shows that $D$ is flat. Note that in general
$D$ is flat if and only if $M$ has constant sectional curvature.

\subsection{The closed 1-form associated to an infinitesimal deformation of 
a hyperbolic metric}

Starting from a $1$-parameter family of hyperbolic metrics on $M$ we have 
constructed above an automorphic vector field $u$ on $\Mt$. Let $\s_u$ be 
its canonical lift. We claim that $\s_u$ is itself automorphic, in the sense that
$\gamma^* \s_u-\s_u$ is a parallel section in $E$. Indeed, since the group of 
deck transformations acts by isometries on $\Mt$,
it commutes with $\mathrm{curl}$, hence from \eqref{canlif}
it also commutes with the linear operator $\s$:
\[\s_{\gamma^* u}=\gamma^*\s_u, \text{ for all } \gamma\in\pi_1(M).\]
We have seen in Section \ref{sec2.1} that
$\kappa:=\gamma^* u-u$ is Killing. Thus by Lemma \ref{lkf}
\[\gamma^*\s_u-\s_u=\s_{\gamma^* u}-\s_u=\s_{\kappa}\]
is parallel as claimed.

Let $\dd$ be the de Rham differential twisted by the flat connection $D$.
Let $\tilde\omega$ be the $1$-form 
\[\tilde\omega=\dd \s_u.\]
Since $\dd$ commutes with the action of $G$, we see that $\tilde\omega$ is 
$G$-invariant on $\Mt$:
\[\gamma^*\tilde\omega=\dd\gamma^*\s_u=\dd(\s_u+\s_\kappa)=\tilde\omega.\]
Thus $\tilde\omega$ descends to a $1$-form $\omega$ on $M$ with values in 
$E=T_\C M$. This form is closed since by construction it is locally exact.

\subsection{Link between infinitesimal deformations 
and $H^1(M,E)$.}

Let us gather below a few facts about $\omega$.

\subsubsection*{The closed $1$-form $\omega$ does not depend on the choice of 
the family of developing maps $\dev_{g_t}$.} 

Indeed, if we replace $\dev_{g_t}$
by $\dev_{g_t}'=a_t\dev_{g_t}$, then $u'=u+\kappa$ where 
$\kappa=\dev_{g_0}^*\dot{a_t}$ is a Killing  field, so 
\[\dd \s_{u'}=\dd \s_u+\dd \s_\kappa=\omega.\]

\subsubsection*{The $1$-form $\omega$ is exact if and only if the 
infinitesimal deformation $\gd$ of the hyperbolic metric is trivial.}

In one direction this is clear:
a vector field $u$ on $M$ determines a germ of a $1$-parameter group
of local diffeomorphisms $\Phi_t$; choose $g_t:=\Phi_t^* g$.
Then $\dev_{g_t}$ may be chosen as 
$\dev_0\circ \tilde{\Phi}$ so the vector 
field of the deformation will be precisely the lift of $u$ to $\Mt$. Thus
$\s_{\tilde{u}}$ is the lift to $\Mt$ of the section $\s_u$ in $E$ over 
$M$ defined by \eqref{canlif}, in other words $\omega$ is exact already on $M$. 
Conversely, assume that there exists $\alpha\in \cun(M,E)$ with 
$\dd\alpha=\omega$. Lifting to $\Mt$ we get
\[\dd\tilde\alpha=\dd \s_u\]
so by Lemma \ref{lkf}, there exists a Killing vector field $\kappa$ on $\Mt$
with $\tilde\alpha-\s_u=\s_\kappa$. Replace $\dev_{g_t}$ by 
$\dev_{g_t}'=a_t\dev_{g_t}$, where $a_t=\exp(t{\dev_g}_*\kappa)$ is a family
of isometries of $H^3$. 
Thus $\tilde\alpha=\s_{u+\kappa}=\s_{u'}$. Since
$\tilde\alpha$ is $G$-invariant, so must be $u'$, therefore $u'$ defines a
vector  
field on $M$, which by definition means that the deformation is trivial.

\subsubsection*{Any closed form $\alpha\in\Lambda^1(M,E)$ is cohomologous 
to $\dd \s_u$ for an automorphic vector field $u$ on $\Mt$} 

Indeed, the lift of $\alpha$ to $\Mt$ is exact, since 
$\Mt$ is simply connected. Thus $\tilde{\alpha}=\dd a$. Now 
decompose $\cun(\Mt,E)$ as follows:
\[\cun(\Mt,E)=\mathrm{Range}(\s)\oplus i\cun(\Mt,T\Mt)\]
where $\s$ is the differential operator \eqref{canlif}. Note that both
spaces in the right-hand side are $G$-invariant. With respect to 
this decomposition we write $a=\s_{u}+iv$, so
\[\tilde{\alpha}=\dd a=\dd \s_{u}+i\dd v.\]
Let $\gamma\in G$; since $\tilde{\alpha}=\dd a$ is invariant, it follows that
$a$ is automorphic, thus there
exists a Killing vector field $\kappa$ with $\gamma^*a-a=\s_\kappa$. 
Put this together with
\begin{align*}
a=\s_{u}+iv,&& \gamma^* a=\s_{\gamma^* u}+i\gamma^* v.
\end{align*}
Since $\mathrm{Range}(s)$ and $i\cun(\Mt,T\Mt)$ are transversal in 
$\cun(\Mt,E)$ we deduce
\begin{align*}
\gamma^* u-u=\kappa,&& \gamma^* v-v=0.
\end{align*}
Hence $iv$ descends to a section of $E$ on $M$; by subtracting $\dd$ 
of this section from $\alpha$ we get the cohomologous form $\dd \s_{u}$
with $u$ automorphic as required.

In summary, we have shown that the application
\[\gd\mapsto [\omega]_{H^1(M,E)}\]
is a well-defined isomorphism between the space of infinitesimal 
deformations of the hyperbolic structures on $M$ modulo trivial deformations, 
and $H^1(M,E)$. Note that the argument holds more generally 
for deformations of metrics of constant sectional curvature of any sign.

\subsection{The variety of representations}
\label{ssc:variety}

To go from infinitesimal deformations (as in Theorem \ref{tm:rigidity})
to small deformations (as in Theorem \ref{tm:existence})  it is
necessary to understand the structure of the space of representations
of $\pi_1(M_r)$ in $PSL(2,\C)$ in the neighborhood of the holonomy representation
$\rho$ of a convex co-compact manifold with particles $M$.

We call $R(M_r)$ the space of representations
of $\pi_1(M_r)$ in $PSL(2,\C)$, that is, the space of homomorphisms of $\pi_1(M_r)$
in $PSL(2,\C)$. The representation $\rho$ is irreducible: note (as in \cite{bromberg2}) that the 
restriction of $\rho$ to each boundary component of $M$ is the holonomy
representation of a complex projective structure on the complement of
the singular points, and as such it is irreducible because any reducible
representation fixes a point in $\C P^1$ and is therefore the holonomy
representation of an affine structure.

It then follows from a result of Thurston \cite{thurston-notes,culler-shalen} 
(see \cite{HK}, Theorem 4.3) that the irreducible component of $R(M_r)$ containing
$\rho$ is a complex variety.

Following a classical construction of Weil \cite{weil:remarks}, one can associate
to $R(N)$ a scheme $\cR(N)$, based on the choice of a presentation of $\pi_1(N)$. Then
the Zariski tangent space of $\cR(N)$ at $\rho$ is naturally associated to the space
of 1-cocycles $Z^1(\pi_1(N); Ad(\rho))$ (see \cite{HK}, Proposition 4.1).

We will see in subsection \ref{ssc:local} that $R(M_r)/PSL(2,\C)$ is actually a smooth complex
manifold in the neighborhood of $\rho$, and that its tangent space is canonically
identified with $H^1(M, E)$.

\section{The geometry of convex co-compact manifolds
with particles} \label{gccccm}

The goal of the next two sections is to find a convenient way to ``normalize''
infinitesimal deformations of convex co-compact manifolds with particles close to
infinity. This will then be used to prove an 
infinitesimal rigidity statement
for hyperbolic manifolds with particles, Theorem
\ref{tm:rigidity}. Here ``normalize'' means to write them 
as sections of a certain bundle which are in $L^2$.  
Theorem \ref{tm:rigidity} will then
follow from an analytical argument; this argument was originally due to Calabi
\cite{Calabi} and Weil \cite{weil}, and has been extended recently to
hyperbolic cone-manifolds by Hodgson and Kerckhoff \cite{HK}. The treatment
here of deformations close to infinity is inspired by the recent work of 
Bromberg \cite{bromberg1,bromberg2}, while the general approach is related 
to the argument used by Weiss \cite{weiss-local}. It would also be interesting to
compare the methods used here to the ones developed by Montcouquiol
\cite{montcouquiol1,montcouquiol-these,montcouquiol2} 
to treat similar questions in higher dimensions, in the
setting of Einstein manifolds with conical singularities.

\subsection{Hyperbolic cone-manifolds}

\subsubsection*{Definitions.}

Hyperbolic cone-manifolds were defined by Thurston \cite{thurston-notes},
using a recursive definition. We define first the special case when 
the singular locus is a disjoint union of lines (i.e.\ a graph without vertices).
Consider a fixed, oriented  
hyperbolic geodesic $\Delta_0$ in $H^3$, and let $U$ be the universal cover
of the complement of $\Delta_0$ in $H^3$. Let $V$ be the metric completion of
$U$, so that $V\setminus U$ is canonically identified with $\Delta_0$; it will
be called the {\it singular set} of $V$. For each
$\alpha>0$, let $V_\alpha$ be the quotient of $V$ by the rotation of angle
$\alpha$ around $\Delta_0$; the image under this quotient of the singular set
of $V$ is called the singular set of $V_\alpha$. Another description of 
$V_\alpha$ is as the hyperbolic cone over the spherical surface $S_\alpha$ with two cone
singularities, both of angle $\alpha$: 
$$ V_\alpha = (S_\alpha\times \R_{>0}, dt^2 + \sinh^2(t) h)~, $$
where $h$ is the metric on $S_\alpha$.

A {\it hyperbolic cone-manifold} with singular locus a union of lines is
a complete metric space for which each point has a neighborhood which is
isometric of an open subset of $V_\alpha$, for some $\alpha>0$. The points
which have a neighborhood isometric to an open subset of the complement of the
singular set in $V_\alpha$ are called {\it regular points}, and the others {\it
  singular points}. The set of regular points of a hyperbolic cone-manifold is
a (non-complete) hyperbolic 3-manifold.

We are interested here in a more general notion of cone-manifolds, for which
the singular set is a graph, that is, three singular lines can meet at a
``vertex''. We require however the number of vertices and edges to be
finite. Such cone-manifolds are made of three kinds of points: in addition to
the regular points and to the points of the singular lines, already
described above, there can be ``vertices'', i.e., 
points which have a neighborhood isometric to a
hyperbolic cone over a 2-dimensional spherical cone-manifold (see
\cite{thurston-notes}). Given such a vertex $v$ in a cone-manifold $M$,
the 2-dimensional spherical cone-manifold over which the neighborhood of 
$v$ in $M$ is ``built'' is called the {\it link} of $v$. 
Each singular point of the link of $v$
corresponds to one of the singular lines ending at $v$, 
and the angle around the singular point in the link 
is equal to the angle around the corresponding singular curve.

Here we suppose that the angle at each singular line is less than $\pi$, it
follows that the same condition holds at each singular point of the link of
$v$, so that the corresponding singular curvature is larger than $\pi$. 
So it follows from the Gauss-Bonnet theorem that the link of each vertex
can have at most 3 cone points. The picture is further simplified by
the fact that a spherical cone-manifold with 3 singular points where the
singular curvature is positive is the double cover of a spherical
triangle (this is a special case of a theorem of Alexandrov, see
\cite{alex,luo-tian}). 

Let $M$ be a hyperbolic cone-manifold. 
Each singular point $x$ of $M$ other than the vertices
has a neighborhood which is isometric to a
subset of $V_\alpha$ for a unique $\alpha>0$; we call $\alpha$ the {\it angle}
of $M$ at $x$. By construction, the angle is locally constant on the singular
lines of $M$.

We will consider here only hyperbolic cone-manifolds which are homeomorphic to
the interior of a compact manifold with boundary, with the singular
set sent by this homeomorphism to an embedded graph, with a finite number
of vertices, and with the exterior vertices 
on the boundary. 

We can follow the definition of $V_\alpha$ above using the Poincar\'e model of
$H^3$, taking as $\Delta_0$ the intersection with the ball of a line $D_0$
going through the origin. This leads to a
conformal model of $V_\alpha$: $V_\alpha$ is conformal to the quotient by a
rotation of angle $\alpha$ of the universal cover of the complement of $D_0$
in the Euclidean ball of radius $1$. 

\subsubsection*{Hyperbolic manifolds with particles}

The specific class of hyperbolic cone-manifolds that we consider contains 
the convex co-compact hyperbolic manifolds, as well as analogous
cone-manifolds, see Definition \ref{df:coco}. Those convex co-compact  
cone-manifolds
are required to contain a non-empty, compact, geodesically convex subset. 

The properties of convex subsets in hyperbolic manifolds with particles are
quite reminiscent of the corresponding properties in non-singular hyperbolic
manifolds. Some considerations on this can be found in the appendix.
One key property, which we will need here, is Lemma \ref{lm:ends}.

Consider
a convex subset $K$ in a convex co-compact cone-manifold $M$, let $N^1(K)$
be the unit normal bundle of $K$, as defined in the Appendix. $N^1K$ contains
all unit vectors in $TM$ which are orthogonal to a support plane of $K$ (and
oriented towards the exterior), as well as some vectors based at the intersection
of $\dr K$ with the singular set of $M$ (see Definition \ref{df:n1}). When
$\dr K$ is smooth and ``orthogonal'' to the singular locus, $N^1K$  is 
homeomorphic to $\dr K$, and diffeomorphic outside the singular points, but
in general $N^1K$ is only a $C^{1,1}$ surface, it has one singular point for 
each intersection point between $\dr K$ and the singular set of $M$. (The
$C^{1,1}$ structure on $N^1K$ is clear if $K$ is the convex core of $M$ and
the support of its bending lamination is a disjoint union of closed curves. 
It is almost as clear if $\dr K$ is polyhedral. In the general case its
existence follows from a limiting argument.)

We consider the restriction of the exponential map to the normal bundle of 
$K$, as the map
$$ \exp_K: N^1K\times (0,\infty)\rightarrow M $$
such that $\exp_K(v,s)=\exp(sv)$, where $\exp$ is the usual exponential map.

The content of Lemma \ref{lm:ends} is that, for any non-empty, convex, compact subset $K$ of
$M$, $\exp_K$ is a homeomorphism from $N^1K\times (0,\infty)$ to $M\setminus K$. 

\subsection{Induced structures at infinity.}

Let $M$ be a hyperbolic manifold with particles, let $M_r$ and
$M_s$ be
the subsets of its regular and of its singular points, respectively.  
$M_r$ has a natural (non-complete) hyperbolic metric, and its universal cover
$\Mt_r$ has a locally isometric projection $\dev$ to $H^3$ which is 
unique up to composition on the left by an isometry of $H^3$. 

The metric completion of $\Mt_r$ is the union of $\Mt_r$ with a union 
of connected sets, each of which
projects to a connected component of $M_s$ and also, by $\dev$, to a complete
graph in $H^3$.

Let $\dr_\infty H^3$ be the boundary at infinity of $H^3$. Then $\dev$ has a
natural extension as a local homeomorphism:
$$ \devt: \Mt_r\cup \dr_\infty \Mt_r \rightarrow H^3\cup\dr_\infty H^3~, $$
where $\dr_\infty \Mt_r$ can be defined, as $\dr_\infty H^3$, as the space of
equivalence classes of geodesic rays in $\Mt_r$, for which the distance to the
singular locus is bounded from below by a positive constant, 
where two rays are in the same class if and only if they are asymptotic. 

The boundary at infinity of $H^3$ can be canonically identified to $\C
P^1$, so that $\devt$
induces on $\dr_\infty \Mt_r$ a complex projective structure. 
We get 
the same $\C P^1$-structure if we compose $\devt$ to the left by an isometry.
Furthermore,
since the hyperbolic isometries act on $\dr_\infty H^3$ by
complex projective transformations, the fundamental group of $M_r$ acts on
$\Mt_r$ by hyperbolic isometries which extend to $\dr_\infty \Mt_r$ as complex
projective transformations. Therefore, $M_r$ has a well-defined boundary at
infinity, which is the quotient of $\dr_\infty \Mt_r$ by the fundamental group
of $M_r$, and which carries a canonical complex projective structure.
 
Let $K\subset M$ be a compact convex subset. The map $\exp:N^1K\times (0,\infty)
\rightarrow M\setminus K$ can be used to define a ``limit'' 
$\exp_\infty:N^1K\rightarrow \dr_\infty M$ (technically, the image of a point 
$(x,v)\in N^1K$ is the equivalence class of the geodesic ray $t\mapsto 
\exp((x,v),t)$). Lemma \ref{lm:ends} shows that this map is a homeomorphism
from $N^1K$ to $\dr_\infty M$, by construction it sends the singular points 
of $N^1K$ to the endpoints at infinity of the singular curves of $M$. 
This shows in particular that two cone singularities in $M$ end at 
different points in $\dr_\infty M$ (i.e., they can not be asymptotic). 
The following statement is a consequence.

\begin{lemma} \label{lm:sing-inf}
Each point $x\in \dr_\infty M$ has a neighborhood which is isometric
either to a half of the hyperbolic space $H^3$ (when $x$ is not an endpoint
of one of the singular curves) or to a neighborhood of
one of the endpoints of $D_0$ in the Poincar\'e model of $V_\alpha$
described above (when $x$ is an endpoint of one of the singular curves,
of angle $\alpha$).
\end{lemma}

\begin{proof}
Suppose that $x$ is the endpoint at infinity of a singular ray $p$ in the
singular set of $M$. Since, by Lemma \ref{lm:ends}, $K$ intersects all
singular rays in $M_s$, $p$ intersects $\dr K$ at a point $y$. Let $n$
be the singular point in $N^1K$ corresponding to the intersection with $p$, 
so that the projection of $n$ to $\dr K$ is $y$ and $n$ is directed along $p$.

Let $U$ be a neighborhood of $n$ in $N^1K$, and let $C=\exp(U\times (0,\infty))$. 
Then $C$ contains a cylinder of exponentially expanding radius around $p$ --
this follows from standard arguments on the normal exponential map of a convex
surface in $H^3$ -- and the statement of the lemma follows.
\end{proof}

Clearly the complex projective structure at infinity is defined on $\dr_\infty M_r$
only, and does not extend to the endpoints at infinity 
of the singular curves in $M$. An extension is however possible. For this
note that $\dr_\infty M_r$ is projectively equivalent, in the neighborhood
of an endpoint at infinity $x$ of a cone singularity, to a neighborhood
in the boundary at infinity of $V_\theta$ of one of the endpoints of the
singular line (here $\theta$ has to be equal to the angle at the singular
line ending at $x$). Considering such model neighborhoods leads to a natural
notion of ``complex projective structure with cone singularities''.

One can also consider the conformal structure underlying the complex projective
structure at infinity ; we will call it the conformal structure at infinity of $M$. It
is defined in the complement, in $\dr N$, of the points which are
the endpoints of the singular graph. (Recall that $N$ is 
the compact manifold with boundary introduced after Definition \ref{df:coco}).
We will see in Remark \ref{removsing}
that this conformal structure can be extended to
the singular points, hence it can also be considered as a 
conformal structure on $\partial N$ with some marked points.

\subsection{The $L^2$ deformations.} \label{ssc:l2}

The regular set $M_r$ of a hyperbolic cone-manifold $M$ carries by definition
a (non-complete) hyperbolic metric. The deformation theory outlined in section 2
for hyperbolic manifolds therefore applies to this setting. There is a 
natural vector bundle over $M_r$, which we still call $E$,
with fiber at a point the vector space of Killing fields
in a neighborhood of this point. Moreover $E$ can be identified 
with $T_\C M_r$ with its natural metric and the
flat connection \eqref{D} (see also \cite{HK}), which we
still call $D$, with flat sections the sections corresponding to a fixed
Killing field. Finally, the infinitesimal deformations of the hyperbolic
cone-manifold structure are associated to closed
1-forms with values in $E$, with two
1-forms corresponding to equivalent deformations if and only if the difference
is the differential of a section of $E$.

Let $\omega$ be a closed 1-form
on $M_r$ with values in $E$. Then $\omega$ is in $L^2$ if:
$$ \int_{M_r} \| \omega\|^2_E dv <\infty~, $$
where the norm of $\omega$ is measured with respect to the hyperbolic
metric on $M$ and the natural metric, at each point of $M$, on $E$.
The tensor product connection $\nabla\otimes D$, where $\nabla$
is the Levi-Civit\`a connection of $M$, can be applied to $\omega$, to obtain a tensor $D\omega$
whose norm can also be measured with respect to the same metrics; again,
$D\omega$ is $L^2$ if the integral of the square of its norm converges over
$M$. 
 
The following lemma is a key point of this paper. It is proved at the end of
the next section, after some preliminary constructions, since it uses some
details on the normalization of a deformation near the singular set and in the
neighborhood of infinity.

\begin{lemma} \label{lm:h1}
Let $(M,g)$ be a hyperbolic manifold with particles. Let $\gd$ be
a infinitesimal deformation of $g$, among hyperbolic manifolds with
particles, which changes neither the conformal
structure at infinity nor the angles at the singular arcs. Then there is a
deformation 1-form $\omega$ associated to $\gd$ which is $L^2$ and such that
$D\omega$ is $L^2$.
\end{lemma}

Here again, the conformal structure at infinity which is considered is the
conformal structure with marked points corresponding to the endpoints of the
singular lines. Note that the fact that $D\omega$ is in $L^2$ is not used in
the sequel, it is included in the lemma since it follows from the proof and
because it could be of interest in different situations.

\subsubsection*{Convex surfaces close to infinity.}

It will be useful, in order to obtain a good normalization of the infinitesimal
deformations of the hyperbolic metrics close to infinity, to find a
foliation of the ends by convex surfaces which are ``orthogonal'' to the
singularities. We first consider another notion of ``convex core'' containing
the singular locus of $M$. We suppose from here on that $M$ is
not one of the model spaces $V_\theta$ defined above.

\begin{defi}
The smallest convex subset of $M$ containing $M_s$ is called 
the \emph{singular convex core} of $M$
and is denoted by $C_S(M)$.
\end{defi}

Here ``smallest'' should be understood for the inclusion; the existence of
$C_S(M)$ is clear since $M$ itself is convex, and the intersection of two
convex subsets of $M$ containing $M_s$ is itself convex 
in the sense of Definition \ref{df:convex},
it cannot be empty since it always contains $M_s$.

Close to infinity, $C_S(M)$ is ``thin'' and concentrated near the singular
locus, as stated in the next proposition. For each $r>0$, we define $CC_r(M)$
as the set of points of $M$ which are at distance at most $r$ from the convex
core $CC(M)$. It is not difficult to check that, for any $r>0$, $CC_r(M)$ is
convex (this follows from the arguments in the appendix).

\begin{prop} \label{pr:expo}
There exists a constant $C>0$ such that, for each $r>0$, any point $x\in
C_S(M)$ which is not in $CC_r(M)$ is at distance at most $Ce^{-r}$ from the
singular locus.
\end{prop}

\begin{proof}
In the Poincar\'e model of $V_\alpha$ described above, 
the intersection of the model with a Euclidean 
ball, with boundary orthogonal to the boundary of the model, 
which does not intersect the singular segment, is isometric to a
hyperbolic half-space. Considering such balls which are tangent to the
singular segment at its endpoint, and which are small enough to fit in the
neighborhood of the endpoint which appears in 
Lemma \ref{lm:sing-inf}, we can find for each endpoint $x_\infty$ 
of the singular graph
$M_s$ in $M$ a finite set of half-spaces $H_1, \ldots, H_p\subset M$,
disjoint from the singular set $M_s$, such that any point $y$ which is at
distance at most $ce^r$ from $M_s\setminus CC_r(M)$
but not in $\cup_{i=1}^n H_i$ is actually at distance at most $c' e^{-r}$ from
$M_s$ (for some $c, c'>0$).

By construction, $C_S(M)$ is contained in the complement of the $H_i, 1\leq
i\leq p$. It follows that, maybe after changing the constants $c$ and $c'$, 
any point in
$C_S(M)\cap (M\setminus CC_r(M))$ which is at distance at most $ce^r$ from one
of the singular curves is actually at distance at most $c'e^{-r}$ from this
singular curve.

However, for $r$ large enough, a point $y\in M\setminus CC_r(M)$ which is at
distance at least $ce^r$ from all the singular curves can not be contained in
$C_S(M)$, since one can construct a half-space in $M$, disjoint from the
singular locus, which contains it. The statement follows.
\end{proof}

\begin{defi}
Let $\Sigma$ be a surface in $M_r$, and let $\Sigmab$ be its closure as a
subset of $M$; suppose that $\Sigmab\setminus \Sigma\subset M_s$. 
We say that $\Sigmab$ is {\it orthogonal to the singular locus}
if, for each $x\in \Sigmab\setminus \Sigma$ and each sequence $(x_n)_{n\in \N}$ 
of points of $\Sigmab$ converging to $x$, the ratio between the distance from $x_n$ 
to $x$ in $\Sigmab$ and the distance from $x_n$ to $M_s$ in $M$ converges to $1$.
\end{defi}

\begin{lemma} \label{lm:foliation}
There exists a constant $k_0>0$ (e.g. $k_0=2$) and 
a compact, convex subset $K\subset M$, with $K\supset CC(M)$, such
that the complement of $K$ in $M$ is foliated by equidistant surfaces, 
which are smooth and locally convex outside the singular locus, 
with principal curvatures at most equal to $k_0$, 
and orthogonal to the singular locus. 
\end{lemma}

\begin{proof}
Choose $r>0$, on which more details will be given below. Let 
$M_s$ be the singular graph
of $M$, and let $x$ be a point in $M_s$ at distance $r$ from
$CC(M)$. By the definition of a cone-manifold given above, there is a
neighborhood $\Omega$ of $x$ in $M$ which is isometric to a ball $\Omega'$
centered at a point
$y$ of $\Delta_0$ in $V_\alpha$, for some $\alpha\in \R_+$. 

Recall that the universal cover of $V_\alpha\setminus \Delta_0$ has a
canonical projection to the complement of a line (which we also call
$\Delta_0$) in $H^3$. The metric completion of the universal cover of
$V_\alpha\setminus \Delta_0$ is obtained by adding a line, which we still call
$\Delta_0$, which contains a unique point $y'$ corresponding to $y$.

Let $y''$ be the image of $y'$ in $H^3$, and let $Q$ be the
plane orthogonal to $\Delta_0$ at $y''$. Then the lift of $Q$ to the universal
cover of $V_\alpha\setminus \Delta_0$ is a totally geodesic subspace $Q'$
which is orthogonal to $\Delta_0$ at $y'$. $Q'$ projects to $V_\alpha$ as
a totally geodesic subset $Q$ which is also orthogonal to $\Delta_0$ at
$y$.

We call $P$ the subset of $\Omega\subset M$ which corresponds to the subset
$Q\cap \Omega'$ of $\Omega'\subset V_\alpha$. If $r$ is large enough,
Proposition \ref{pr:expo} indicates that 
$C_S(M)\setminus P$ has two connected components, one of which is contained in an
$\epsilon$-neighborhood of the
subset of $\gamma$ which is bounded by $x$ on the side opposite to $CC(M)$.

Since the same construction can be done for each of the points at distance $r$
from $CC(M)$ in the singular locus of $M$, we can ``cut out'' the neighborhoods
in $C_S(M)$ of the parts of the singular curves which are at distance more
than $r+1$ from $CC(M)$. Since this is done by cutting along totally geodesic
surfaces which are orthogonal to the singular locus, we obtain in this way a
compact subset $K'$ of $M$, contained in $C_S(M)$ and in $C_{r+1}(M)$, which
is convex. However the boundary of $K'$ is not smooth.

We can now call $K$ the set of points of $M$ at distance at most $1$ from
$K'$; it is again compact and convex, and its boundary is $C^{1,1}$
smooth and strictly convex. 
Smoothing this boundary surface by any of the classical techniques
-- without changing it in a neighborhood of its intersections with the
singular curves, where it is totally umbilic -- yields a convex, compact
subset $K$ of $M$ with a boundary which is smooth and orthogonal to the
singular locus. The statement is then obtained by considering the foliation of
the complement of $K$ by the surfaces at constant distance from $K$.

Consider an integral curve of the unit vector field orthogonal to these
surfaces, towards infinity. Since the surfaces are equidistant, this integral
curve is a geodesic, and a classical computation (see e.g. \cite{GHL}) shows
that, along it, the second fundamental form of the surfaces satisfies a
Riccati equation: 
$$ B' = I - B^2~. $$
It follows that the principal curvatures of the equidistant surfaces converge
to $1$ close to infinity in each of the ends of $M$. Therefore, replacing $K$
by a larger compact subset if necessary, we obtain that the principal
curvatures of the equidistant surfaces are at most $k_0$, for any choice of
$k_0>1$.  
\end{proof}

In the sequel, for each end $e$ of $M$, we call $(S_{e,t})_{t\in \R_+}$ the
family of surfaces obtained in the previous lemma, which foliates a
neighborhood of infinity in the end $e$.

\subsubsection*{The metric at infinity associated to an equidistant foliation.}

Such an equidistant foliation $(S_{e,t})_{t\in \R_+}$ determines a 
natural metric $g_{\infty, e}$ on the connected component of $\dr_\infty M$
corresponding to $e$, it is defined as:
$$ g_{\infty,e} = \lim_{t\rightarrow \infty} e^{-2t} I_t~, $$
where $I_t$ is the induced metric on $S_{e,t}$. The surfaces $(S_{e,t})_{t\in
  \R_+}$, and the boundary at infinity, are identified through the orthogonal
projections on the $S_{e,t}$. The homothety factor $e^{-2t}$ is designed to
compensate the divergence of $I_t$ as $t\rightarrow \infty$.

Clearly, the conformal structure of
$g_{\infty,e}$ is equal to the conformal structure underlying the $\C
P^1$-structure on $\dr_\infty M$ which was already mentioned above. It also
follows quite directly from its definition that $g_{\infty,e}$ is a smooth
metric with conical singularities at the endpoints of the singular lines of
$M$, where its singular angle is equal to the singular angles around the
corresponding singular lines of $M$. 

Note that $g_{\infty, e}$ is not in general hyperbolic, it depends on the
choice of the equidistant foliation $(S_{e,t})_{t\in \R_+}$. Actually it is
possible to choose $g_{\infty, e}$ and deduce from it an equidistant
foliation, which might however only be defined for $t\geq t_0$, for some
$t_0\in \R$ (see e.g. \cite{fefferman-graham,graham-lee} where related
questions are treated in the more general context of conformally compact
Einstein manifolds, but without singularities, or \cite{volume} for the
3-dimensional hyperbolic case).

It is perhaps worth noting that there is another possible definition of
the metric at infinity $g_{e,\infty}$: it is equal to $e^{-2t} I^*_t$, where
$I^*_t$ is the ``horospherical metric'' of $S_{e,t}$,
i.e., $I_t+2\II_t+\III_t$ (where $I_t, \II_t$ and $\III_t$ are the induced
metric, second and third fundamental forms of $S_{e,t}$, respectively) for any
choice of $t$ --- the result does not depend on the choice of $t$. 
Details on this can be found in \cite{c-epstein,horo}.

\subsubsection*{Geodesics close to infinity.}

A direct consequence of the existence of the foliation by parallel, convex
surfaces orthogonal to the singular locus, obtained in the previous paragraph,
is the existence of another foliation, by geodesics going to infinity and
normal to those surfaces. 

\begin{lemma}
For each end $e$ of $M$, for each $x\in S_{e,0}$, there exists a geodesic ray
$h_{e,x}$ with endpoint $x$ which is orthogonal to the surfaces $S_{e,t}, t\in
\R_+$. The geodesic rays $h_{e,x}, x\in S_{e,0}$, foliate $\cup_{t\in \R_+}
S_{e,t}$, and, for each $x\in S_{e,0}$, the point at distance $t$ from $x$ in
$h_{e,x}$, called $h_{e,x}(t)$, is in $S_{e,t}$. 
\end{lemma}

\begin{proof}
This is a direct consequence of the previous lemma, taking as the $h_{e,x}$
the curve orthogonal to the equidistant surfaces and starting from $x$.
\end{proof}

\begin{lemma} \label{lm:geodesics}
For each end $e$ of $M$, there exists a constant $C_e>0$ with the following
property. Let $\gamma:[0,1]\rightarrow S_{e,0}$ be a smooth curve, then, for
each $t\in \R_+$:
$$ \left\| \frac{\dr h_{e,\gamma(s)}(t)}{\dr s}\right\| \geq C_e \|
\gamma'(s)\| e^t~. $$
\end{lemma}

\begin{proof}
The $h_{e,\gamma(s)}$ are geodesics, and are orthogonal to $\gamma$. Moreover:
$$ \frac{\dr}{\dr t} 
\left\| \frac{\dr h_{e,\gamma(s)}(t)}{\dr s}\right\|_{|t=0} \geq 0~, $$
because $S_{e,0}$ is convex. So the
estimate follows directly from classically known estimates on the behavior of
Jacobi fields along a geodesic, see e.g., \cite{GHL}.
\end{proof}

\section{The normalization of infinitesimal deformations}

The goal of this section is to prove Lemma \ref{lm:h1}. The argument uses some
additional notations, which we first introduce. We denote by $\Mt$ the universal
cover of $M$, with its singular locus. So $\Mt$ is a quotient of the metric
completion of the universal
cover $\Mt_r$ of $M_r$, already defined above, with ramification at the lift
to $\Mt$ of the singular locus $M_s$. The boundary at infinity $\dr_\infty \Mt$
of $\Mt$ can
be defined in the same way as the boundary at infinity of $\Mt_r$ 
(as the space of geodesic
rays up to the equivalence relation ``being asymptotic''). The complement
of the singular points in $\dr_\infty \Mt$ is the
quotient of the complement of the endpoints of the singular curves in the 
boundary at infinity of $\Mt_r$ by the group acting on $\Mt_r$ 
with quotient the complement of the singular curves in $\Mt$. 

There are three main steps in the proof of Lemma \ref{lm:h1}.
The first is to normalize a family of hyperbolic cone-metrics
$g_s$ with cone angles constant in $s$ by a family of isotopies, so that the
automorphic vector field $v$ on $\Mt$ associated to the deformation extends
to an automorphic vector field $V$ on $\dr_\infty \Mt$. Moreover, $V$ will
turn out to be equivalent to a holomorphic vector field $V+W'$, where $W'$ is the
lift to $\dr_\infty \Mt$ of a vector field defined on $\dr_\infty M$,
and the behavior 
of $W'$ near the singular points of $\dr_\infty M$ can be understood
thoroughly. 

The second step is to construct from $V+W'$ a section $F$ of a bundle of 
quadratic polynomials on $\dr_\infty \Mt$, which is strongly 
related to the bundle $E$ of local Killing fields on $\Mt$, and
use the description of $W'$ at the singular points to show that $F$
also behaves rather nicely close to the singular points. 

Finally the third step uses the section $F$ to construct a deformation
1-form $\omega$ in $M$ equivalent to the initial deformation. 
The estimates on $F$ then translate as the required estimates on $\omega$.

\subsection{The vector field at infinity}

\begin{lemma}\label{lnorm}
Let $g_s$ be a $1$-parameter family of hyperbolic cone-metrics on $M$ with 
constant angles at the singular graph. There 
exists a $1$-parameter family $\Phi_s$ of isotopies of $M$ such that for all 
$s$, 
\begin{enumerate}
\item the hyperbolic cone-metric $g_s':=\Phi_s^*g_s$ coincides
with $g_0$ in a model neighborhood of the singular locus near 
infinity, and near the vertices;
\item the geodesic half-lines $h_{e,x}(t)$, defined using $g_0$, are also
geodesic for $g_s'$.
\item the deformation $1$-form for the family $g_s'$ is uniformly 
bounded in norm near the singular graph, and vanishes 
in a model neighborhood of the singular half-lines outside the convex core. 
\end{enumerate}
\end{lemma}
\begin{proof}
From Lemma \ref{lm:sing-inf}, if the cone angles are fixed then
the metrics $g_s$ are all isometric in a neighborhood of the singular 
lines near infinity to a subset in a fixed model neighborhood $V_\alpha$.
Also near the 
vertices of the singular graphs, hyperbolic metrics with fixed cone 
angles are rigid (see Proposition \ref{locd}). However the metrics $g_s$ may
vary in a neighborhood of the singular graph, even if the 
angles are fixed: there may appear an elongation of the singular segments, 
and also a twist of the graph along such a segment.

If we fix these lengths and the twists of the graph along segments, the 
hyperbolic cone-metric is clearly
rigid near the singular graph. We construct now 
some explicit metrics $g_s'$ with the same lengths and twists as $g_s$, which 
are therefore isometric to $g_s$ on 
the $\epsilon$-neighborhood $U_\epsilon$
of the singular graph. 
Let $l_{g_s}(e)$ be the length of the 
edge $e$ with respect to $g_s$, and $\theta_{g_s}(e)$ the additional 
twisting angle of $g_s$ along $e$, as compared to $g_0$. 
Choose $\epsilon$ sufficiently small
so that the singular graph is a deformation-retract
of its $\epsilon$-neighborhood $U_\epsilon$. Cut this neighborhood 
into pieces using totally geodesic disks orthogonal to the singular graph
at some fixed distance $\delta$ from the vertices. We obtain in this way
for each edge $e$ a finite-length cone-manifold $C_e$ of angle $\alpha(e)$
around a singular curve of length $l_{g_0}(e)-2\delta$. Replace this cylinder
by the cylinder of the same type of length $l_{g_s}(e)-2\delta$, and glue
it back with a twist of angle $\theta_s$.

One can realize this metric on $U_\epsilon$ (outside the singular locus) 
as follows: fix an edge $e$, 
let $l_s$ be the length of $e$
with respect to $g_s$, also $\theta_s$ the twist of $g_s$ along $e$ 
(relative to $g_0$). Let $(x,z)$
be coordinates adapted to $g_0$ on the cylinder corresponding to $e$, 
where $x\in [\delta,l_0-\delta]$ is the height function
and $z\in \C$ is a complex variable in the disk of radius $\epsilon$, written
$z=(r,\theta)$ in polar coordinates:
\[g_0=dr^2+\cosh(r)^2dx^2 +\sinh(r)^2d\theta^2.\]
Let $\phi$ be a cut-off function on 
$[0,l_0]$, which vanishes for $x<\delta$, is increasing, 
and equals $1$ for $x>l_0-\delta$. 
Pull back (at time $s$) the metric $g_0$ through the map
\begin{equation}\label{cd}
(x,r,\theta)\mapsto (x+\phi(x)(l_s-l_0), r, \theta+\phi(x)\theta_s).
\end{equation}
These maps for different edges do not glue nicely to $U_\epsilon$
(they do not agree near the vertices) but the pulled-back metric
does extend to $U_\epsilon$, and is isometric to $g_s$ since they have the same 
elongation and twist along each edge. 

Since the metrics $g_s$ and $g_s'$ are isometric on $U_\epsilon$, we 
can pull-back $g_s$ through a family of isotopies of $M$ 
starting from the identity at $s=0$,
such that 
the resulting metrics are equal to $g_s'$ near the singular locus. The surface
$S_{e,0}$ constructed in Section \ref{gccccm}
is convex also for the metrics $g_s'$ for sufficiently small $s$.
Choose a second family of isotopies which is the identity on a neighborhood 
of the singular locus near infinity and on the convex core, and which 
maps the normal geodesics flow from $S_{e,0}$ (with respect to $g_s$) 
onto the corresponding flow with respect to $g_0$. The points where the 
surfaces $S_{e,0}$ (for different values of $s$) intersect the singular lines
may lay at varying distance from the convex core; we choose this second 
family of isotopies to be given by \eqref{cd}
between these intersection points and the convex core, with $\theta_s=0$
and the necessary elongation $l_s$.

The metrics $g_s'$ coincide near the singular locus at infinity and near the 
vertices, hence the associated deformation vector field $v$ is
 Killing so the deformation $1$-form $\omega$ vanishes in the above region 
as claimed. It remains to check that $\omega$ is bounded on the cylinders 
near each singular 
segment $e$. Recall that $\omega=\dd \s_v$, where $\s$ is the canonical
lift operator \eqref{canlif}, and $v$ is the vector field 
tangent in $s=0$ to the $1$-parameter family of maps \eqref{cd}:
\[v=\phi(x)(\dot{l}\partial_x+\dot{\theta}\partial_\theta)\]
where the constants $\dot{l}$, $\dot{\theta}$ are the infinitesimal 
variations of the length and twist of the edge $e$. It is
not hard to see
that $\partial_x$ and $\partial_\theta$ are Killing fields, which correspond
to translations along $e$, respectively rotations around $e$, therefore
$\dd \s_{\partial_x}=\dd \s_{\partial_\theta}=0$. Also one sees easily
that 
\begin{align*}
\s_{\phi\partial_x}=\phi\s_{\partial_x},&&\s_{\phi\partial_\theta}=
\phi\s_{\partial_\theta}-\frac{i}{2}\phi'(x)\tanh(r) \partial_r.
\end{align*}
It follows that
\[\omega=\phi'(x) dx\otimes(\dot{l}\s_{\partial_x}
+\dot{\theta}\s_{\partial_\theta})-\frac{i\dot{\theta}}{2}
\dd(\phi'(x)\tanh(r)\partial_r).\]
The volume form is $\sinh(r)\cosh(r)dr d\theta dx$. A straightforward computation 
shows that $\omega$ is bounded uniformly near $e$.
\end{proof}

\begin{remark}\label{db}
The same computation shows that $D\dd\s_v$  is also uniformly bounded at 
finite distance from the convex core.
\end{remark}

\subsubsection*{The deformation vector field in the ends.}

Let $v$ be the automorphic vector field on $\tilde{M}_r$ defined as in section
\ref{local} from a family of hyperbolic cone-metrics with constant cone
angles, normalized as in Lemma \ref{lnorm}. Note that the normalization 
from Lemma \ref{lnorm} gives in particular an identification of the 
boundaries at infinity for the different metrics. 

Let $E$ be
an end of $M$, i.e., a connected component of the complement in $M$ of a 
non-empty, compact, convex subset $K$. The singular set of $E$ is a
disjoint union of singular rays $p_1, \cdots, p_N$. The boundary at infinity 
$\dr_\infty \Mt_r$ is the
disjoint union of (possibly countably) copies of the $\dr_\infty E_r$, where  
$E$ runs through all the ends of $M$. Consider $v$ on 
each such copy of $\Et_r$; as can be checked locally in hyperbolic space
$H^3$, it extends smoothly to $\dr_\infty \Et_r$.
Moreover $v$ is Killing in a neighborhood of 
the singular locus near infinity, since the family $g_s$ was 
normalized to be constant there. 
We call $V$ the automorphic vector field on $\dr_\infty \Mt_r$ obtained by
extending $v$ in this manner. Since $v$ is Killing in a neighborhood of 
the singular lines near infinity, it follows that
$V$ is locally a projective vector field near the singular points.

\subsubsection*{A holomorphic vector field on $\dr_\infty \Mt$.}

We will need an elementary and well-known statement: given an automorphic 
vector field on a Riemann surface, it is equivalent to a holomorphic automorphic
vector field if and only if the induced infinitesimal variation of the complex
structure (considered up to isotopy) vanishes. 

\begin{lemma} \label{lm:complex}
Let $\Sigma$ be a closed surface 
with marked points $x_1,\ldots, x_n$, endowed with a $\C P^1$-structure $\sigma$ 
with singularities at the $x_i$. Set 
$\Sigma_r:=\Sigma\setminus\{x_1,\ldots,x_n\}$.
Let $\phi:\tilde\Sigma_r\to\C P^1$ be the developing map 
of $\sigma$, and let $V$, a section of $T\tilde\Sigma_r$, be an automorphic
vector field corresponding to a infinitesimal variation of $\sigma$ (among the
$\C P^1$-structures). Suppose that the infinitesimal variation of the complex
structure on $\Sigma$, marked by the position of the $x_i$,
vanishes up to isotopy. Moreover, suppose that $V$ is projective in the lift 
of a uniform 
neighborhood of the singular points. Then there exists a
smooth vector field $W$ on $\Sigma$ (i.e., smooth at the $x_i$), vanishing 
at the $x_i$, such that if $W'$ is the lift of $W$ to $\Sigmat_r$, 
then $V+W'$ is a holomorphic vector field.
\end{lemma}

The smoothness of $W$ at $x_i$ is to be understood for the underlying
complex structure on $\Sigma$. Since $V$ is projective, in particular
it does not change, at first order, the angle around the singular points.

\begin{proof}
Let $J$ be the complex structure underlying the $\C P^1$-structure $\sigma$.
By our hypothesis, $V$ does not change the complex structure --- marked by the
position of the singular points --- on $\Sigma$, considered up 
to diffeomorphisms isotopic to the identity. This means precisely
that the action of $V$ on the complex structure is the same as the action of a
vector field defined on $\Sigma$, which we call $-W$, which vanishes at
the singular points. Calling $W'$ the lift of $W$ to $\tilde{\Sigma}$, it 
is clear that $V+W'$ does not change pointwise the complex structure on 
$\Sigma$ (again, marked by the position of the singular points)
so that $V+W'$ is a holomorphic vector field. It follows that 
$W$ is holomorphic in the neighborhood where $V$ is projective. 
\end{proof}

Since $v$ was normalized to be Killing near infinity in a neighborhood
of the singular locus, it follows that
$V$ is indeed projective near the singular points.
It follows from the previous lemma that we can replace the vector field
$V$ on $\dr_\infty \Mt$ by another vector field $V+W'$, corresponding to 
the same infinitesimal variation of the $\C P^1$-structure, but which is 
holomorphic.

\subsection{A vector bundle of quadratic polynomials}

We recall here some well-known notions on a natural bundle of polynomials of
degree at most 2 on a surface with a complex projective structure.

\subsubsection*{Complex polynomials and Killing fields.}

It is necessary to understand the relationship (partly based on the
Poincar\'e half-space model) between hyperbolic Killing fields, projective
vector fields on $\C P^1$, and polynomials of degree at most
2 over $\C$ (or in other terms holomorphic vector fields over
$\C P^1$). 

\begin{remark}
Let $\kappa$ be a Killing field on $H^3$. Let $\kappab$ be the image of
$\kappa$ in the Poincar\'e half-space model. Then $\kappab$ has a continuous
extension as a vector field on the closed half-space $\{ z\geq 0\} $. On the
boundary $\{ z=0\}$, the restriction of this extension is tangential to the
boundary plane, and its coordinates are given --- after identification of $\{
z=0\}$ with $\C$ --- by a polynomial of degree at most 2.
\end{remark}

\begin{proof}
Let $(\Phi_t)_{t\in [0,1]}$ be a one-parameter family of hyperbolic
isometries, with $\Phi_0=I$. For all $t\in [0,1]$, let $\phi_t$ be the action
of $\Phi_t$ on the boundary at infinity, identified with
$\overline{\C}$. Then, for all $t\in [0,1]$, $\phi_t$ acts on $\C$ as:
$$ \phi_t(z) = \frac{a(t)z+b(t)}{c(t)z + d(t)}~, $$
with $a(0)=d(0)=1, b(0)=c(0)=0$. Taking the derivative at $t=0$, we find that:
$$ \left(\frac{\dr}{\dr t} \phi_t(z)\right)_{|t=0} = (a'(0)z + b'(0)) - z
(c'(0)z + d'(0))~, $$
and the result follows.
\end{proof}

In other words, 
the hyperbolic Killing fields act on the boundary at infinity of
$H^3$, identified with $\C P^1$, as holomorphic vector fields. 
Moreover, given any point $z_0\in \C P^1$, $\C P^1\setminus \{ z_0\}$
can be identified with $\C$, and can therefore be given a complex
coordinate $z$. The action at infinity of the Killing fields are of the
form: 
$$ v(z) = P(z) \dr_z~, $$
where $P$ is a polynomial of degree at most $2$. The set of 
these polynomials is
invariant under the action of the M\"obius transformations, so that the notion
of polynomial of degree at most $2$ makes sense on any surface endowed with a
$\C P^1$-structure. 
More details on the relation between quadratic polynomials and Killing vector
fields can be found in \cite{bromberg1,bromberg2}.

\subsubsection*{Estimates on Killing fields in terms of polynomials.}

The different monomials have a simple interpretation in terms of hyperbolic
Killing fields:
\begin{itemize}
\item Polynomials of degree 0 correspond to Killing fields that vanish at 
  the point at infinity in $\overline{\C}$, and fix (globally) 
  the horospheres ``centered'' at this point at infinity.
\item Homogeneous polynomials of degree 1 correspond to Killing fields that
  fix (globally) the hyperbolic geodesic corresponding, in the Poincar\'e
  half-space model, to the vertical line containing $0$. They are sums of
  infinitesimal rotations around this geodesic and infinitesimal translations
  along it.
\item Homogeneous polynomials of degree 2 correspond to Killing fields that
  vanish at the origin, and fix (globally) the horospheres ``centered'' at
  this point.
\end{itemize}

These three types of Killing fields, and their interpretation, have a direct
generalization to the more general situation of a hyperbolic 3-manifold $M$,
in terms of the behavior at infinity, near a point $z_0\in \dr_\infty \Mt$, of
the Killing vector fields defined on $\Mt$. We consider an affine complex
coordinate $z$ defined in the neighborhood of $z_0$, i.e., the actions at
infinity of the Killing vector fields are of the form $P(z) \dr_z$, where $P$
is a polynomial of degree at most $2$. 

\begin{lemma} \label{lm:est-poly}
There exists a constant $C>0$ with the following property. Let $x\in H^3$ and
let $N\in T_xH^3$ be a unit vector such that $\lim_{s\rightarrow \infty}
\exp_x(sN)=z_0\in \dr_\infty H^3$. Let $P$ be the totally geodesic plane
orthogonal to $N$ at $x$, let $g_0$ be the induced metric on $P$, and let
$G:P\rightarrow \dr_\infty H^3$ be the hyperbolic Gauss map. Suppose that,
at $z_0$, $G_* g_0 = e^{2r} |dz|^2$. Then: 
\begin{itemize}
\item the Killing vector field $\kappa_1$ corresponding to the polynomial
  $(z-z_0)\dr_z$, considered as a flat section of $E$, has norm bounded, at
  $x$, by $C$. 
\item the Killing vector field $\kappa_2$ corresponding to the polynomial
  $(z-z_0)^2\dr_z$ has norm bounded, at $x$, by $Ce^{-r}$.
\end{itemize}
\end{lemma}

The norm which is considered here is not the norm of Killing fields,
considered as vector fields on $H^3$, but rather their norm considered as
(flat) sections of the vector bundle $E$; recall that this norm depends on the
point of $H^3$ where they are considered.

\begin{proof}
Both statements follow from a direct computation, for instance using the
Poincar\'e half-space model.
\end{proof}

Clearly the previous statement could be extended to include Killing vector
fields corresponding to polynomials of degree $0$, however this will not be of
any use here. It is also worth noting that a possible proof uses the
invariance under the multiplication of $z-z_0$ by a constant $\lambda$; then
$(z-z_0)\dr_z$ does not change, while $(z-z_0)^2\dr_z$ is multiplied by
$\lambda$. Under the same homothety, $\kappa_1$ does not change along the
``vertical'' geodesic ending at $z_0$, while
$\kappa_2$ is multiplied by $\lambda$ because it corresponds to a parabolic
isometry fixing the horospheres ``centered'' at $z_0$.

\subsubsection*{The vector bundle of quadratic polynomials.}

The remarks in the previous paragraph lead naturally to define a bundle over
$\C P^1$, which is strongly related to the bundle of local Killing fields,
which is used on $H^3$ or on any hyperbolic manifolds. Although the definition
is given here on $\C P^1$, it should be clear that it is of a local nature,
and makes sense for any surface with a $\C P^1$-structure.

\begin{defi}
We call $P$ the trivial bundle over $\C P^1$,  with fiber at each point the
vector space of holomorphic vector fields on $\C P^1$. 
\end{defi}

Clearly $P$ has a natural flat connection $D^P$, such that the flat sections
are those which correspond, at each point of $\C P^1$,
to the same holomorphic vector field. In other terms, $D^P$ is the trivial connection
on the trivial bundle $P$.

\subsubsection*{The section of $P$ associated to a vector field.}

Given a vector field on $\C P^1$, or more generally on a surface with a $\C
P^1$-structure, one can associate to it a section of the bundle $P$, defined
by taking at each point the ``best approximation'' by polynomial vector fields
of degree at most $2$.

\begin{defi} \label{df:section}
Given a holomorphic vector field $v$ defined on an open subset
$\Omega\subset \C P^1$, there is a
section $F$ of $P$ which is naturally associated to $v$; at each point $z_0\in
\Omega$, $F_{z_0}$ is equal to the holomorphic vector field on $\C P^1$ which best
approximates $f$. Given any affine identification of $\C P^1$ (minus a point) 
with $\C$, if $v:=f\dr_z$, this translates in $\C$ as:
$$ F_{z_0}(z) = (f(z_0) + (z-z_0) f'(z_0) + \frac{(z-z_0)^2}2
f''(z_0))\dr_z~. $$ 
\end{defi}

\begin{lemma} \label{lm:vanishing}
Let $v$ be a holomorphic vector field on $\Omega\subset \C P^1$. 
Let $F$ be the associated
section of $P$. Then, at each point $z_0\in \Omega$, $D^PF$ has values in the
subspace of $P_{z_0}$ of vector fields which vanish, along with their first
derivatives, at $z_0$. 
\end{lemma}

\begin{proof}
Since the statement is local, the proof
takes place in $\C$, and we write $v=f(z)\dr_z$. 
Let $z_0\in \Omega$ , and let 
$Z\in T_{z_0}\C$. We identify vector 
fields on $\C$ with complex functions on $\C$ and obtain, 
using the definition of the flat connection
$D^P$, that for all $z$ in some open subset of $\C$:
\begin{align*}
  (D_Z^PF)(z_0) = & (Z\dr_{z_0} F_{z_0}(z))\dr_z \\
  = & Z\dr_{z_0}\left(f(z_0) + (z-z_0) f'(z_0) +
  \frac{(z-z_0)^2f''(z_0)}{2}\right)\dr_z \\
  = &Z \left(f'(z_0) - f'(z_0) + (z-z_0) f''(z_0) - (z-z_0)f''(z_0) +
  \frac{(z-z_0)^2f'''(z_0)}{2}\right)\dr_z, 
\end{align*}
so that:
\begin{equation}
  \label{eq:der1}
(D_Z^PF)(z_0) = \left(Z \frac{(z-z_0)^2f'''(z_0)}{2}\right)\dr_z~. 
\end{equation}
This shows that $(D^PF)(z_0)$ 
takes its values in the vector space of homogeneous polynomials of degree $2$,
as needed.
\end{proof}

\subsection{The geometry of $\dr_\infty M$ near the singular points.}

We now concentrate on an explicit description of the complex structure and
complex projective structure on $\dr_\infty M$ near its singular points, which
will be necessary in estimates below.

\subsubsection*{The boundary at infinity of $\Mt$.}

We have already noted that the boundary at infinity of $\Mt$ carries a $\C
P^1$ structure, with singular points corresponding to the endpoints of the
singular arcs. It also carries a vector bundle, $P$, with fiber at each
point $x$ the vector space of vector fields in the neighborhood of $x$ which
are obtained as continuous extensions to the boundary (for instance in a local
Poincar\'e model) of hyperbolic Killing vector fields. 

By the (local)
considerations above, the fiber of $P$ at $x$ can also be identified with the
vector space of projective vector fields in a neighborhood of $x$. Again, $P$
has a natural flat connection, still called $D^P$, 
with flat sections the sections corresponding 
to a given projective vector field. Since its statement is of a local nature,
Lemma \ref{lm:vanishing} still holds on $\dr_\infty M_r$.

\subsubsection*{Special coordinates near the singular points.}

We now consider more carefully what happens on the boundary at infinity 
of $M$ in the neighborhood of a singular point. Let $x_1,\cdots, x_n$ be
the singular points on $\dr_\infty M$, i.e., the endpoints of the 
singular arcs. For each $i\in \{ 1,\cdots, n\}$, the $\C P^1$-structure
of $\dr_\infty M_r$ in the neighborhood of $x_i$ is projectively 
equivalent to a neighborhood of the vertex in a ``complex cone'' 
which we call 
$\C_{\theta_i}$: it is the quotient of the universal cover of
the complement of $0$ in $\C$ by a rotation of center $0$ and angle 
$\theta_i$, where $\theta_i$ is the angle at $x_i$. 

We choose a neighborhood $\Omega_i$ of $x_i$, and a complex projective map
$u:\Omega_i\rightarrow \C_{\theta_i}$ sending $x_i$ to $0$, which is a
diffeomorphism from $\Omega_i\setminus \{ x_i\}$ to its image. The map $u$ is
uniquely determined (by the complex projective structure) up to composition on
$\C_{\theta_i}$ with a rotation and a homothety; we choose this homothety so
that, as $x\rightarrow x_i$, the metric $g_\infty$ on $T_x\dr_\infty M$
behaves as $|du|^2$. It follows that there is a constant $C$, independent of
$i$, such that, on the $\Omega_i\setminus \{ x_i\}$:
$$ \frac{|du|^2}{C} \leq g_\infty \leq C|du|^2~. $$

There is a natural holomorphic local 
diffeomorphism from $\C_{\theta_i}$ to $\C$. With
obvious notations, it is defined by sending a point $u\in \C_{\theta_i}$ to
$u^{2\pi/\theta_i}$. With the same notations we set $z:=u^{2\pi/\theta_i}$,
this defines a complex coordinate $z$ on $\Omega_i$. 
Therefore we have proved:

\begin{remark}\label{removsing}
The boundary at infinity 
$\dr_\infty M$ can be canonically considered as a smooth surface, 
with a smooth complex structure; only the complex projective structure 
and the metric at infinity, $g_\infty$, have singularities at the endpoints 
of the singular arcs of $M$. (The ``singularities'' of the complex projective
structure are in the sense explained just before subsection \ref{ssc:l2}.) 
\end{remark}

\subsubsection*{Estimates on the deformation field at infinity.}

We now have most of the tools necessary to 
``normalize'' the infinitesimal deformations
of the hyperbolic structure of a manifold with particles. This
means that, given a infinitesimal deformation $\gd$ of the metric $g$ keeping
the cone angles and the conformal structure at infinity fixed, 
we will show that it is
associated to a 1-form $\omega$ with values in the bundle $E$ which is of a
very special form. It will then follow that 
$\omega$ and $D\omega$ are in $L^2$. 

The first step is to associate to $\gd$ an automorphic vector field $v$, along
the ideas at the end of section 2, using the infinitesimal deformation of the
development map. We have seen that $v$ can be chosen to have a continuous
extension to the boundary at infinity of the universal cover of $M$. 
We call $V$ the automorphic vector field on $\dr_\infty \Mt$ obtained by
extending $v$ in this manner. Since $v$ can be chosen to be Killing near the 
singular curves, it follows that $V$ is projective near the singular points
of $\dr_\infty \Mt$.

According to Lemma \ref{lm:complex}, there exists a holomorphic vector field 
$V+W'$ on $\dr_\infty \Mt_r$ such that $W'$ is the lift to 
$\dr_\infty \Mt_r$ of a vector field $W$ defined on $\dr_\infty M$.
It is then clear that $V+W'$ is automorphic. Note that $W'$ is holomorphic in a 
neighborhood of the singular points and vanishes at the singular points. Indeed,
$V$ is projective, hence holomorphic, in a small enough neighborhood, therefore
$W'$ itself must be holomorphic. Moreover, by construction 
$W$ preserves the \emph{marked} smooth structure of $\dr_\infty M$, which means 
that it vanishes at the singular points.

Choose $i\in \{ 1,\cdots, n\}$.
Let $F$ be the section of $P$, associated by Definition \ref{df:section} to
$V+W'$. We use the coordinate $u$ on the $\Omega_i$ defined above, so a vector
tangent to $\dr_\infty M$ can be identified with a complex number.

\begin{lemma} \label{lm:47}
Let $i\in \{ 1, \cdots, n\}$ and choose $u_0\in \Omega_i\setminus \{ x_i\}$. 
For all vectors fields $U, U'$ defined in a neighbourhood of 
$u_0\in \dr_\infty\Mt$ we have: 
$$ D_U^PF (u_0) = U\alpha(u_0) (u-u_0)^2 \dr_u~, $$
while
$$ D_{U'}^PD_U^P F (u_0) = UU'(\beta(u_0) (u-u_0)^2 + \gamma(u_0) (u-u_0))\dr_u~. $$
Moreover, $\alpha$ and $\gamma$ are bounded by a constant $C>0$, and there
exists another constant $\epsilon_0\in (0,1)$ such that:
\begin{equation} \label{eq:beta} 
|\beta(u_0)| \leq C/|u_0|^{1-\epsilon_0}~. 
\end{equation}
\end{lemma}

\begin{proof}
On compact sets disjoint from the singular points, the estimates follow directly 
from Lemma \ref{lm:vanishing}. Thus we consider only points
$u_0$ in a neighborhood of $x_i$ where $V$ is projective. 
Now $W'$ is holomorphic near $x_i$ and vanishes at $x_i$, so it admits
a Taylor series decomposition
\[W' = (w_1 z + w_2 z^2 + w_3 z^3 \cdots)\dr_z \]
in a neighborhood of $x_i$.

Let $u$ be an ``affine coordinate'' at $x_i$ for the $\C P^1$-structure 
induced on $\dr_\infty M$, as defined above.
Let $\mu:=2\pi/\theta$, then $\mu>2$ since $\theta\in (0,\pi)$. 
Then $dz=\mu u^{\mu-1} du$, so that $\dr_z = \mu^{-1} u^{1-\mu} \dr_u$, and
it follows that:
\begin{equation*}
W' = (w_1 u + w_2 u^{1+\mu} + w_3 u^{1+2\mu} + \cdots) \dr_u~.
\end{equation*}
The section
of $P$ associated to $V$ is parallel by definition on the set where $V$ is 
projective, so it is enough
to estimate the covariant derivatives of the section associated
to the vector field $W'$. Thus we may assume that $F$ is the section in $P$ 
associated to $W'$.
Equation (\ref{eq:der1}) shows that:
$$ D^P_UF_{u_0} = \frac{U}{2}\left(w_2 \mu(\mu^2-1)u_0^{\mu-2}
+ 2w_3\mu(4\mu^2-1)u_0^{2\mu-2} + \cdots\right)(u-u_0)^2 \dr_u~. $$
Taking one more differential leads to:
$$
D^P_{U'}D^P_U F_{u_0} = - UU'\left(w_2 \mu(\mu^2-1)u_0^{\mu-2}
+ 2w_3\mu(4\mu^2-1)u_0^{2\mu-2} + \cdots\right)(u-u_0) \dr_u $$
$$ +
\frac{UU'}{2}\left(w_2
  \mu(\mu^2-1)(\mu-2)u_0^{\mu-3} + 2w_3\mu(4\mu^2-1)(2\mu-2) u_0^{2\mu-3} + 
  \cdots\right)(u-u_0)^2 \dr_u~.
$$
Moreover, $\mu>2$, so the estimates announced in
the lemma follow directly from the two previous
equations, by taking $\epsilon_0:=\mu-2$.
\end{proof}

It follows by compactness from this statement that the same estimates hold 
on $\dr_\infty M$, with Eq.\ (\ref{eq:beta}) replaced by $|\beta(u_0)| \leq C$
in the complement of the union of the $\Omega_i, 1\leq i\leq n$, with $u$
taken to be an ``affine'' coordinate, compatible with the complex projective
structure, defined on a finite number of compact domains covering $\dr_\infty
M$. 

\subsection{The deformation 1-form in the ends}

It remains now to define from the automorphic section $F$ a section of 
$E$ over $\Mt_r$ which is also
well-mannered close to infinity, in the same way as in \cite{bromberg1}.

\subsubsection*{Normalization from infinity.}

By Lemma \ref{lm:geodesics}, 
there is a compact subset $K\subset M$ whose complement
is foliated by geodesic rays $h_{e, x}$, with $x\in \dr K$. These rays lift 
to geodesic rays $\hti_{e,x}$ in $\Mt_r$, where $e$ is a lift of an end of $M_r$
and $x\in \dr\Kt_r$.

We now define a
section $\kappa$ of $E$ over $\Mt_r\setminus \Kt$ as follows: for each $y\in
\Mt_r\setminus \Kt$, let $e,x$ be the unique elements such that $y\in
\hti_{e,x}$. Then $\kappa_x$ is the Killing field (defined in the neighborhood of
$h_{e,x}$) with extension at infinity (in the neighborhood of the endpoint $z$
of $\hti_{e,x}$) the projective vector field corresponding to $F_z$. 
Clearly the section $\kappa$ of $E$ defined in this way over 
$\Mt_r\setminus \Kt$ is smooth. Let $\s_v$ be the canonical lift of the 
deformation vector field $v$ to a section of $E$.
By lemma \ref{lm:complex}, $\s_v-\kappa$ is $G$-invariant on the ends. 
Let $\phi(t)$ be a cut-off function depending
on the distance function $t$ to the convex core, which vanishes for $t\leq 1$
and equals $1$ for large $t$. Then $\phi(t)(\s_v-\kappa)$ is well-defined and
$G$-invariant on $\Mt_r$. Thus $\s_v-\phi(t)(\s_v-\kappa)$ is automorphic, 
differs from $\s_v$ by a $G$-invariant section and behaves near infinity 
like $\kappa$. Consider the invariant 1-form 
$\omega:=\dd(\s_v-\phi(t)(\s_v-\kappa))$ on $\Mt_r$ with values in $E$. 
By construction, this form and the initial $1$-form
$\dd(\s_v)$ correspond to equivalent infinitesimal deformations of the hyperbolic 
cone-manifold structure on $M$. Moreover, both
$\omega$ and $D\omega$ vanish in the direction of the lines $h_{e,x}$ near
infinity.

\subsubsection*{Different metrics on $\dr_\infty M$.}

It is natural to consider, on the boundary at infinity of $M$, the
metric $g_\infty$ which was already defined above in terms of the foliation 
of the ends near infinity. On the leafs of this foliation, however, there
are two metrics which are quite natural:
\begin{itemize}
\item the ``horospherical metric'' $I^*_t:=I_t+2\II_t+\III_t$. It is conformal
  --- through the Gauss map --- to the metric $g_\infty$ at infinity,
\item the metric $g_t$ which is defined as follows. For each 
  $x\in S_{e,t}$, let $P_x$ be the totally geodesic plane tangent to
  $S_{e,t}$ at $x$, then $T_xS_{e,t}=T_xP_x$, and the metric $g_t$, on
  $T_xS_{e,t}$, is equal to the pull-back of $g_\infty$ to $P_x$ through the
  Gauss map $G:P_x\rightarrow \dr_\infty M$.
\end{itemize}
Note that $g_t$ is not equal to the pull-back to $S_{e,t}$ of $g_\infty$ by
the Gauss map $G:S_{e,t}\rightarrow \dr_\infty M$. However each is bounded
by a constant times the other. Recall that $k_0$ was 
defined above as an upper bound on the principal curvatures of the 
surfaces $S_{e,t}$.

\begin{remark} \label{rk:comp-metriques}
For all $t\in \R_+$, we have:
\begin{enumerate}
\item for all $x\in S_{e,t}$, if $G:S_{e,t}\rightarrow \dr_\infty M$ is the
hyperbolic Gauss map, then:
$$ g_t \leq G^*g_\infty \leq (1+k_0)^2 g_t~, $$
\item $I^*_t = e^{2t}G^*g_\infty$, where $G:S_{e,t}\rightarrow \dr_\infty M$
  is the hyperbolic Gauss map.
\end{enumerate}
\end{remark}

\begin{proof}
In the first point, the first inequality follows from the convexity of
$S_{e,t}$, because the differential on $T_xS_{e,t}$ of the Gauss map of 
$S_{e,t}$ is ``larger'' than the differential of the Gauss map of the
totally geodesic plane tangent to $S_{e,t}$ at $x$. The second inequality
follows in the same way from the fact that the principal curvatures of
$S_{e,t}$ are bounded by $k_0$.

The second point is a direct consequence of the fact, already mentioned above,
that the horospherical metric changes in a very simple way along an
equidistant foliation (see \cite{horo}).
\end{proof}

Moreover, the metric $g_t$ is the one appearing in Lemma \ref{lm:est-poly},
which yields an estimate in terms of $t$ of the Killing vector fields
associated to special quadratic polynomials on $\dr_\infty M$.

\begin{cor} \label{cr:est}
There exists a constant $C'>0$ as follows.
Let $x\in S_{e,t}$, and let $u_0:=G(x)$, where $G$ is the hyperbolic Gauss map
of $S_{e,t}$. Let $u$ be an affine coordinate system in the  neighborhood of
$u_0$ (for the $\C P^1$-structure on $\dr_\infty M)$, chosen so that
$|du|^2=g_\infty$ on $T_{u_0}\dr_\infty M$. Then the Killing vector field
$\kappa_1$ corresponding to the vector field $(u-u_0)\dr_u$ (considered as a
flat section of $E$) has norm, at $x$,
bounded by $C'$, while the Killing vector field $\kappa_2$ corresponding to
the vector field $(u-u_0)^2\dr_u$ has norm bounded, at $x$, by $C'e^{-t}$. 
\end{cor}

\begin{proof}
Direct consequence of Lemma \ref{lm:est-poly}.
\end{proof}

\subsubsection*{Estimates on $\omega$ and $D\omega$.}

It is now possible to estimate the $L^2$ norm of $\omega$, and then of
$D\omega$, so as to prove Lemma \ref{lm:h1}. Let $x\in S_{e,t}$, and let $X\in
T_xM$. We are interested in $\omega(X)$, and we already know that $\omega$
vanishes along the lines orthogonal to the surfaces $S_{e,t}$, so we suppose
that $X\in T_xS_{e,t}$. 

Let $U:=G_*X$, where $G:S_{e,t}\rightarrow \dr_\infty M$ is the Gauss map. 
Remark \ref{rk:comp-metriques} shows that: 
$$ \| U\|_{g_\infty}\leq c e^{-t} \| X\|_M~, $$
for some constant $c>0$.

Recall that $\omega(X)=D_X\kappa$, where $\kappa$ is the section of $E$
corresponding to $F$. So $\omega(X)$ corresponds to the vector field $D^P_UF$
on $\dr_\infty M$. According to Lemma \ref{lm:47},
$D^P_UF=U\alpha(u_0)(u-u_0)^2\dr_u$, where $\alpha$ is bounded and $u$ is an
affine coordinate system near $u_0$. Using Corollary \ref{cr:est} we see
that $\omega(X)$ has norm (at $x$) bounded by
$C'' e^{-t} \| U\|_{g_\infty}$, or in other terms by $C'' e^{-2t} \| X\|$,
where $C''>0$ is some constant.

This means that, at $x$,  $\|\omega\|\leq C''e^{-2t}$, so that by Fubini:
$$ \int_0^\infty \int_{S_{e,t}} \|\omega\|^2 da dt \leq 
\int_0^\infty (C'')^2 e^{-4t}
A(S_{e,t}) dt \leq C_3 \int_0^\infty e^{-2t} dt \leq C_3/2~, $$
where $C_3>0$ is yet another constant, and $A(S_{e,t})$ denotes the area of
the surface $S_{e,t}$ (for the ambient metric). Since the same estimate
applies to each end of $M$, we conclude that
$\omega\in L^2$ near infinity.

A similar argument can be used to estimate $D\omega$. Let $X,X'\in
T_xS_{e,t}$, and let $U:=G_*X, U':=G_*X'$. 
As above we have: 
$$ \| U\|_{g_\infty} \leq ce^{-t}\|X\|, ~ \|U'\|_{g_\infty} \leq
ce^{-t}\| X'\|~. $$
But $(D_X\omega)(X')$ corresponds to the vector field 
$$ D^P_UD^P_{U'}F -D^P_{G_*(\nabla_XX')}F $$
on $\dr_\infty M$, which itself can be written using Lemma \ref{lm:47} as:
$$ UU'(\beta(u_0) (u-u_0)^2 + \gamma(u_0) (u-u_0))\dr_u~,$$
where $\gamma$ is bounded and where $\beta(u_0)$ is bounded by $C/|u_0|$ in
the $\Omega_i$ and by $C$ in their complement, 
$C$ being a constant.

This can be written, using Corollary \ref{cr:est}, as the following estimates
when $u_0$ is in one of the $\Omega_i, 1\leq i\leq n$:
\begin{eqnarray*}
\| (D_X\omega)(X')\|_E & \leq & \| U\|_{g_\infty} \| U'\|_{g_\infty} 
\left(\frac{C_4 e^{-t}}{|u_0|^{1-\epsilon_0}} 
+ C_5\right) \\
& \leq & \| X\| \| X'\| e^{-2t} \left(\frac{C_4
    e^{-t}}{e^{-t}d_{S_{e,t}}(x,S_{e,t}\cap M_s)^{1-\epsilon_0}} 
+ C_5\right)  
\end{eqnarray*}
which translates as:
\begin{equation} \label{eq:estimation}
\| D\omega\| \leq e^{-2t} \left(\frac{C_4}{d_{S_{e,t}}(x,S_{e,t}\cap
      M_s)^{1-\epsilon_0}} + C_5\right)~. 
\end{equation}
The same estimates can be used when $u_0$ is not in one of the $\Omega_i$ and
      yields:
\begin{eqnarray*}
\| (D_X\omega)(X')\|_E & \leq & \| U\|_{g_\infty} \| U'\|_{g_\infty} 
\left(C_4 e^{-t} 
+ C_5\right) \\
& \leq & \| X\| \| X'\| e^{-2t} \left(C_4
    e^{-t} + C_5\right)~.  
\end{eqnarray*}
and a simple compactness argument then shows that, perhaps after changing the
constants $C_4$ and $C_5$, Eq.\ (\ref{eq:estimation}) holds on all $S_{e,t}$.

We can now integrate the square of this norm over the ends of $M$, and obtain 
that: 
$$
\int_0^\infty \int_{S_{e,t}} \| D\omega\|^2 dx dt \leq   
\int_0^\infty \int_{S_{e,t}} e^{-4t} \left(\frac{C_4}{d_{S_{e,t}}(x,S_{e,t}\cap
      M_s)^{1-\epsilon_0}} + C_5\right)^2 dx dt~.
$$
Using the comparison between the induced metric $I_t$ on $S_{e,t}$ and the
metric at infinity $g_\infty$, the previous equation translates as:
$$
\int_0^\infty \int_{S_{e,t}} \| D\omega\|^2 dx dt \leq   
\int_0^\infty \int_{\dr_\infty M} e^{-2t} 
\left(\frac{C'_4}{d_{g_\infty}(x,S_\infty)^{1-\epsilon_0}} + 
C'_5\right)^2 dx dt~,
$$
for some constants $C'_4, C'_5>0$, so that, calling 
$S_\infty=\{ x_1,\cdots, x_n\}$ the set of singular points on $\dr_\infty M$:
$$
\int_0^\infty \int_{S_{e,t}} \| D\omega\|^2 dx dt \leq   
\int_0^\infty e^{-2t} dt\int_{\dr_\infty M} \left(
\frac{(C'_4)^2}{d_{g_\infty}(x,S_\infty)^{2-2\epsilon_0}} + 
\frac{2C'_4C'_5}{d_{g_\infty}(x,S_\infty)^{1-\epsilon_0}} + 
(C'_5)^2 \right)dx~.
$$
Note that the area element of $g_\infty$, close to the singular points, 
behaves as $\rho d\rho d\theta$ (where $\rho$ is again the distance to the 
singular points considered). So $1/\rho^{2-2\epsilon_0}$ is integrable, and
it follows that 
all the terms in the integral over $\dr_\infty M$ converge, 
with the contribution of the terms in
$1/d_{g_\infty}(x,S_\infty)^{2-2\epsilon_0}$ and 
$1/d_{g_\infty}(x,S_\infty)^{1-\epsilon_0}$ bounded for each 
singular point of $S_\infty$. (Note that the hypothesis that the angles around
the singular lines are less than $\pi$ is used here.) 
This shows that, for each end of $M$, the
integral of $\| D\omega\|^2$ is bounded. Since this holds for all ends of $M$,
the integral of $\| D\omega\|^2$ is bounded near infinity. 

Now by Lemma \ref{lnorm} and Remark \ref{db}, both $\dd\s_v$ and $D\dd\s_v$
are bounded (hence $L^2$) near the convex core. We still need to check the
integrability of $\phi'(t)dt\otimes(\s_v-\kappa)$
and similarly for the $D$-covariant derivative of this form. The support of 
$\phi'$ is compact and contained in the ends; thus it is enough to check 
that $\s_v-\kappa$ and $D(\s_v-\kappa)$ are bounded in norm on the ends. Recall
that in a neighborhood of a singular line, $\s_v-\kappa$ is the $G$-invariant 
section in $E$ corresponding to the holomorphic vector field $W$ on the boundary 
at infinity. Therefore the required estimate follows again (as in the
arguments right above) from the behavior of the section $F$ on $\dr_\infty
\Mt$ as described in Lemma \ref{lm:47} and from the relation between 
projective vector fields at infinity and Killing fields in $\Mt$ which can be
read from Corollary \ref{cr:est}.

This finishes the proof of Lemma \ref{lm:h1}.

\section{Infinitesimal rigidity}\label{Sir}

In this section we first prove a general result about $L^2$ cohomology
and then we show how to apply it in our setting. 

\subsection{A general argument}

Let $E\to M$
be a vector bundle over a Riemannian manifold together with a flat 
connection $D$ and a Riemannian metric along the fibers. Let 
\[\dd:\Lambda^*(M,E)\to\Lambda^{*}(M,E)\]
denote the twisted de Rham differential with coefficients in $E$, and $\ded$
its formal adjoint. 
Consider the symmetric operator $P:=\dd+\ded: \Lambda^*(M,E)\to\Lambda^*(M,E)$.
Let $L^2$ denote the Hilbert space of square-integrable sections in 
$\Lambda^*(M,E)$. We view $P$ as an unbounded operator with domain $\cunc$, 
the space of smooth compactly supported $E$-valued forms on $M$, which
is dense in $L^2$.

The elements of $L^2$ act as distributions on $\cunc$ and thus they can be
differentiated. For $k\in\Z$ define the Sobolev space $H^k$ as the 
space of those 
section $\phi\in L^2$ such that $P^k\phi\in L^2$ in the sense of distributions 
(or equivalently, $H^k=\Dom (P^k)^*$). 
This is a Hilbert space with the graph norm
\[\|\phi\|^2_{H^k}:=\|\phi\|^2+\|(P^*)^k\phi\|^2.\]
Define also $\hmin^k$ as the completion
of $\cunc$ with respect to the norm $\|\cdot\|_{H^k}$. It is a small lemma that
$\hmin^k$ injects naturally in $H^k$.

By the Friedrichs extension theorem, the operator
\[P^*\oP:\hmin^1\cap H^2\to L^2\] 
is self-adjoint and non-negative. Note that a form belongs to
$\hmin^1\cap H^2$ if and only if its components in all degrees do.

\begin{lemma}\label{lema4.1}
Let $\alpha\in L^2(M,\Lambda^k M\otimes E)$. Assume that
the inequality
\begin{equation}\label{fw}
(\dd+\ded)^2-1\geq 0
\end{equation}
holds on $\cunc(M,\Lambda^k\otimes E)$. Then there exists a unique $\gamma\in
\hmin^1\cap H^2$ of degree $k$ such that
$P^*\oP\gamma=\alpha$. Moreover, if $\alpha$ is smooth then $\gamma$ 
is also smooth.
\end{lemma}
\begin{proof}
Let $(P^*\oP)_k$ denote the restriction of $P^*\oP$ to $k$-forms.
By continuity, Eq.\ \eqref{fw} implies that $(P^*\oP)_k\geq 1$, therefore
$0$ does not belong to its spectrum. In other words, $(P^*\oP)_k$ is 
invertible from the $k$-form part of $\hmin^1\cap H^2$ to $L^2$. 
Finally, if $\alpha$ is smooth then
$\gamma$ is also smooth by elliptic regularity.
\end{proof}

Of course, the lemma holds for any strictly positive constant instead of $1$.

We make now the assumption that $\hmin^1=H^1$, which is another way of saying
that $P=\dd+\ded$ is essentially self-adjoint. Let $\dmin,\dmax,\delmin,
\delmax$ denote the minimal, respectively 
maximal extensions of the operators $\dd$ and $\ded$. An $L^2$ form $\alpha$
is called \emph{closed} if $\dmax\alpha=0$, and \emph{exact} if there exists
$\gamma\in\Dom(\dmin)$ with $\dmin\gamma=\alpha$. 

\begin{prop}\label{cohoegze}
Assume that 
\begin{itemize}
\item the inequality \eqref{fw} holds on $\cunc(M,\Lambda^k M\otimes E)$ 
for $k=0,\ldots,\dim(M)$;
\item the operator $P$ is essentially self-adjoint.
\end{itemize}
Then for each $k$, every closed form $\alpha$ in $L^2(M,\Lambda^k M\otimes E)$ 
must be exact.
Moreover, there exists a primitive $\gamma$ which is in $\hmin^1$, coexact, 
and which is smooth if $\alpha$ is smooth.
\end{prop}
\begin{proof}
The second hypothesis, which translates to $\Dom(P^*)=\hmin^1$, will be used 
throughout the proof without further explanation.

We first claim that $L^2(M,\Lambda^k M\otimes E)$ decomposes orthogonally
into
\[L^2(M,\Lambda^k M\otimes E)=\Ran(\dmin^{k-1})\oplus \Ran(\delmin^{k+1}).\]
Let $\alpha\in L^2(M,\Lambda^k M\otimes E)$. 
By Lemma \ref{lema4.1} and the first hypothesis,
there exists $\beta\in
\hmin^1\cap H^2(M,\Lambda^k M\otimes E)$ such that $P^*\oP\beta=\alpha$.
Clearly $\Dom(\oP)\subset\Dom(\dmin)\cap\Dom(\delmin)$, so $\dmin\beta,\delmin
\beta$ are both well-defined. Furthermore, $\dmin\beta+\delmin\beta$
belongs to $H^1=\hmin^1$, which implies that both $\dmin\beta$ and $\delmin
\beta$ belong to $\hmin^1$ and so $P^*\dmin\beta=\delmin\dmin\beta$,
$P^*\delmin\beta=\dmin\delmin\beta$ are both $L^2$. Hence
\[\alpha=\dmin(\delmin\beta)+\delmin(\dmin\beta)\]
where $\delmin\beta,\dmin\beta$ are both in $\hmin^1$, which proves the claim.

Next, we claim that $\ker\dmax=\Ran(\dmin)$. The inclusion
$\Ran(\dmin)\subset\ker \dmax$ is clear from the definitions.
Let $\alpha=\dmin(\delmin\beta)+\delmin(\dmin\beta)$ be such that 
$\dmax\alpha=0$. Since $\dmin(\delmin\beta)$ is already in $\ker \dmax$, 
we deduce that $\dmax\delmin\dmin\beta=0$.
Also $\Ran(\delmin)\subset\ker\delmax$, so $\delmax\delmin\dmin\beta=0$.
Thus, $\delmin\dmin\beta$ belongs to $\Dom(\dmax)\cap\Dom(\delmax)\subset\Dom
(P^*)=\hmin^1$ and $P^*(\delmin\dmin\beta)=0$. The kernel of $P^*=\oP$
vanishes since $P^*\oP$ is invertible (by the first hypothesis), 
therefore $\delmin\dmin\beta=0$
and so $\alpha=\dmin(\delmin\beta)$ is contained in the range of $\dmin$.

The primitive $\delmin\beta$ is $L^2$ by construction. Moreover, it belongs to
$\Dom(\dmin)$ and to $\ker\delmax$. Thus it belongs to
$\Dom(\dmax)\cap\Dom(\delmax)\subset\Dom(P^*)=\hmin^1$ and 
$P^*(\delmin\beta)=\alpha$. If $\alpha$ is smooth
then by elliptic regularity, $\delmin\beta$ is also smooth.
\end{proof}

We actually proved in particular that the $L^2$ cohomology of 
$M$ twisted by $E$ vanishes.

\subsection{Application to cone-manifolds}
Consider now the bundle $E\simeq TM\oplus TM\simeq T_\C M$ of 
infinitesimal Killing vector fields
on the conical hyperbolic $3$-manifold $M$.
Let $i:E\to E$ be the complex structure, $i(u,v):=(-v,u)$. Define an
endomorphism-valued $1$-form $T:TM\otimes E\to E$ by
\[T_V\phi:=V\times i\phi\]
where $\times$ denotes the vector product in $TM$, acting on each component of 
$\phi$. The flat connexion $D$ on $E$ is given by the explicit formula
\[D_V(u,v)=(\nabla_V u+V\times v, \nabla_V v - V\times u)\]
where $\nabla$ is the Levi-Civit\`a connexion on $TM$. We write this
as $D\phi=(\nabla+T)\phi$.
We extend $T$ and $\nabla$ to $\Lambda^*(M,E)$ as the identity, resp.\ the 
de Rham differential on the form factor. We endow $E$ with the direct sum 
Riemannian metric, and $\Lambda^*M$ with its usual metric.
 
\begin{prop}[Matsushima \& Murakami \cite{MM}]\label{mamu}
The Laplacian of the twisted de Rham differential $\dd$ satisfies 
Eq.\ \eqref{fw} on $\Lambda^k_c(M,E)$ for $k=0,\ldots,3$.
\end{prop}
\begin{proof}
By definition, 
\[\Delta=(\nabla^*\nabla+\nabla\nabla^*)+(T^*T+TT^*)
+T\nabla^*+\nabla^*T+T^*\nabla+\nabla T^*.\]
It was observed by Matsushima \& Murakami that
$T\nabla^*+\nabla^*T+T^*\nabla+\nabla T^*=0$.
Indeed, let $(x_j)_{1\leq j\leq 3}$ be geodesic normal coordinates at a point
$x\in M$, $(e_j)$ the coordinate vector fields and
$(e^j)$ the dual basis. 
Let $\Phi=\alpha\otimes\phi\in\Lambda^*(M,E)$. We have:
\begin{align*}
\nabla\Phi=&d\alpha\otimes \phi+\sum e^j\wedge\alpha\otimes \nabla_{e_j}\phi&
T\Phi=&\sum e^j\wedge\alpha\otimes e_j\times i\phi\\
\nabla^*\Phi=&\delta\alpha\otimes \phi-\sum e^k\lrcorner\alpha\otimes
\nabla_{e_k}\phi&
T^*\Phi=&\sum e^k\lrcorner\alpha\otimes e^k\times i\phi
\end{align*}
where the contraction uses the metric on forms. Since
$\times$ and $i$ commute with $\nabla$, we get at $x$, 
\begin{align*}
(\nabla T^*+T^*\nabla)\Phi=&\sum L_{e_j}\alpha\otimes e_j\times i\phi
+\alpha\otimes e_j\times i\nabla_{e_j}\phi\\
=&-(T\nabla^*+\nabla^*T)\Phi.
\end{align*}
The Laplacian $\nabla^*\nabla+\nabla\nabla^*$ is non-negative.
We claim that $T^*T+TT^*\geq 1$ pointwise. We work at $x\in M$ 
where the basis $e_j$ is orthonormal. 

Let us first examine the action of $T^*T$ on $0$-forms.
It is immediate that 
\[T^*T\phi=\sum_{k=1}^3 e_k\times i(e_k\times i\phi)=2\phi.\]

We focus now on $1$-forms.
Notice that $T^*T$ and $TT^*$
act diagonally on $E$ with respect to the
splitting $E=TM\oplus TM$, so it is enough to prove the claim
on a real section 
\[\Phi=\sum a_{ki} e^k\otimes e_i.\]
Set $\Phi_{il}=e^i\otimes e_l$. Then
\begin{align*}
(T^*T+TT^*)\Phi_{il}=&-\sum_{j,k} e^k\lrcorner( e^j\wedge e^i)
\otimes e_k\times(e_j\times e_l)
+e^j\wedge(e^k\lrcorner e_i)\otimes e_j\times(e_k\times e_l)\\
=:&{\sum}'\Phi_{il}+{\sum}''\Phi_{il}
\end{align*}
where the two sums group the terms with $j=k$, resp.\
$j\neq k$. Since
$e^j\wedge e_j\lrcorner+e_j\lrcorner e^j\wedge=1$, $i^2=-1$
and $\sum e_j\times(e_j\times \phi)=-2\phi$, we find that
$\sum'\Phi_{il}=2\Phi_{il}$.

For $j\neq k$ notice 
that $e^j\wedge e^k\lrcorner+e^k\lrcorner e^j\wedge=0$. Therefore
\begin{align*}
{\sum}''\Phi_{il}=&\sum_{j\neq i}e^j\otimes(e_i\times(e_j\times e_l)
-e_j\times(e_i\times e_l))\\
=&\begin{cases}
\sum_{j\neq i}\Phi_{jj}&\text{if $i=l$};\\
-\Phi_{li}&\text{if $i\neq l$}.
\end{cases}
\end{align*}
In conclusion, for $\Phi=\sum a_{ki} \Phi_{ki}$ we obtain
\[\langle (T^*T+TT^*)\Phi,\Phi\rangle=|\Phi|^2+ \left(\sum_i a_{ii}\right)^2
+\sum_{i\neq k}(a_{ik}-a_{ki})^2\geq |\Phi|^2.\]
Note that the equality is obtained precisely for traceless, symmetric $\Phi$.

For $k=2,3$ we remark that the Hodge $*$ operator commutes with the Laplacian, 
and acts isometrically from $\Lambda^k M\otimes E$ to $\Lambda^{3-k} 
M\otimes E$. Thus the result follows from what we proved above for $k=1,0$.
\end{proof}

\subsection{Cone angles and essential self-adjointness}
The aim of this subsection is to prove that when the cone angles of our 
hyperbolic cone manifold are smaller than $\pi$, the third hypothesis of 
Proposition \ref{cohoegze} is satisfied. The proof is based 
on the analysis from \cite{weiss-local}.

\begin{theorem}\label{esa}
Let $M$ be a hyperbolic manifold with particles. Assume that all cone angles 
belong to the interval $[0,\pi]$. Then the twisted Hodge-de Rham operator 
$P=\dd+\ded$ acting in $L^2(M,\Lambda^*M\otimes E)$ is essentially
self-adjoint. 
\end{theorem}

\begin{proof}
We must show that if $u\in L^2$ and $Pu\in L^2$ then $u\in \hmin^1$.

We first localize $u$ near the singular locus. 
Let $\psi_1:M\to[0,1]$ be a smooth function which equals $1$ 
near the singular graph, and which vanishes 
outside the $\epsilon$-neighborhood of the singular set ($\epsilon$ is chosen 
sufficiently small so that this is a tubular neighborhood). We also need 
$|d\psi_1|$ to be uniformly bounded; 
actually we can choose $\psi_1$ such that
$|d\psi|\to 0$ at infinity, by asking that $\psi_1$ only depends on the distance
function $r$ to the singular set. Then 
$\psi_1 u$ is clearly in $L^2$; moreover, 
\begin{equation}\label{clifford}
P(\psi_1 u)=\psi_1Pu+c(d\psi_1)u
\end{equation}
 is also in $L^2$, where $c$ denotes Clifford multiplication, i.e.,
\[c(\alpha)u=\alpha\wedge u-\alpha\lrcorner u.\] 
Set $\psi_2:=1-\psi_1$.
\begin{lemma}\label{lemaloc}
The form $\psi_2 u$ belongs to $\hmin^1$.
\end{lemma}
\begin{proof}
We rely on the results of Weiss \cite{weiss-local}. 
We follow partly the proof of the fact that on a complete manifold, 
all ``geometric'' differential operators are self-adjoint. For $n\to\infty$ let $f_n$ be a 
smooth function on $M$, equal to $1$ on the $n$-neighborhood of the convex
core, 
supported on the $2n$-neighborhood of the convex core, and such that 
$|df_n|\leq 2/n$. We can choose such a function to depend only on the variable 
$t$ which parametrizes the families of equidistant surfaces on the ends
from Lemma
\ref{lm:foliation}. Clearly $f_n \psi_2 u$ converges in $L^2$ to $\psi_2 u$ as
$n\to\infty$. 
By \eqref{clifford},  
\[P(f_n\psi_2 u)=f_nP(\psi_2 u)+c(df_n)\psi_2u.\]
By Lebesgue dominated convergence, the first term converges
in $L^2$ to $P(\psi_2 u)$ while the second converges to $0$; thus 
$\psi_2 u$ can be approximated by compactly-supported forms in $H^1$ sense, 
as claimed.
\end{proof}

It is left to prove that $\psi_1 u$ 
belongs to $\hmin^1$.
Without loss of generality we can therefore assume that 
$u$ lives in a tubular neighborhood $U_\epsilon$ of the singular graph.

Again by a partition of unity using the cut-off function $f_1$, it is enough to
prove the statement separately for $u$ supported at finite
distance from the singular line in a model conical set $V_\alpha$, and for
$u$ supported near the convex core. 

In the first case, 
we use the Poincar\'e ball model
of $H^3$. Let $g_n$ be a sequence of cut-off functions on the unit
interval (i.e., $g_n:[0,1)\to [0,1]$ is smooth, equals $1$ near $R=0$
and has compact support). Denote by $R$ the radial function on the 
disk. Then $g_n(R)$ is rotation-invariant, so it descends to a function
on $V_\alpha$. We can choose $g_n$ to converge to $1$ on each compact set;
moreover, since the metric $dR^2/(1-R^2)^2$ on $(-1,1)$ is complete,
we can impose that 
\begin{equation}\label{mic}
|dg_n(R)|\leq 1/n.\end{equation}
From \eqref{clifford} we see 
that $g_nu$ and $P(g_nu)$ are both in $L^2$, in other words $g_n u$ belongs to 
$H^1$, the \emph{maximal} domain of $P$. Now $g_n u$ has support inside a ball 
(depending on $g_n$). 

From the results of \cite[Sections 4 and 5]{weiss-local}, 
we claim that $g_n u$ must be in the \emph{minimal} domain of $P$, provided
that 
$\alpha$ is smaller than $\pi$. Indeed, Weiss shows in 
\cite[Proposition 5.10]{weiss-local}
that the bundle $E$ with its connection is \emph{cone-admissible}; this is 
a technical condition which implies \cite[Corollary 4.34]{weiss-local}
that on a compact hyperbolic cone-3-manifold the operator $P$ 
is essentially self-adjoint. Finally, the proof of this last
Corollary is local in nature, and amounts to proving, after the use 
of a partition of unity, exactly the above claim.

As in the proof of Lemma \ref{lemaloc}, we have obtained a sequence
$(g_n u)$ in $\hmin^1$ which converges to $u$ in $L^2$ and which is 
Cauchy in the $H^1$ norm by \eqref{mic}; thus $u$ is itself in $\hmin^1$.

We use the same argument for the remaining case, namely where $u\in H^1$ 
is supported near the convex core. By the local result
\cite[Corollary 4.34]{weiss-local}, again $u$ must belong to $\hmin^1$, provided the
angles are all bounded above by $\pi$.
\end{proof}

\section{Proofs of the main results}

\subsection{Proof of Theorem \ref{tm:rigidity}}
Let $\dot{g}$ be a infinitesimal deformation of a hyperbolic 
metric $g$ with particles among metrics of the same type, which
fixes both the cone angles and the conformal structure at infinity. 
Let $\omega$ be the closed $E$-valued deformation $1$-form 
associated to $\dot{g}$ by Lemma \ref{lm:h1}. We thus know that 
$\omega$ is square-integrable.
The first hypothesis of Proposition \ref{cohoegze} holds by 
Proposition \ref{mamu}. The 
second hypothesis is fulfilled by Theorem \ref{esa} if all
cone angles are at most $\pi$. Thus, by
Proposition \ref{cohoegze}, $\omega$ is exact as a smooth form. 
By the results of Section \ref{local}, the infinitesimal deformation
$\dot{g}$ is trivial.

\subsection{Cohomological arguments}

A curious phenomenon is that sometimes, uniqueness implies existence.
Something similar happens here as we explain below. The arguments used here
are somewhat similar to those in \cite{HK}.

Let $V\subset M_r$ be a compact manifold with boundary which is a deformation
retract of $M_r$. $V$ can be obtain e.g., by smoothing the boundary of the 
complement, inside the convex core, 
of the $\epsilon$-neighborhood of the singular locus.
Let $U$ denote the closure inside $M_r$ of the complement of $V$, in particular 
$U$ is an incomplete manifold with boundary. Note that the natural inclusion map 
on the level of forms induces isomorphisms
\[H_c^k(M_r)\cong H^k(M_r,U).\]
Note that all cohomology groups in this section are twisted by the flat bundle 
$E$, unless otherwise specified; we suppressed $E$ from the notation.
Consider the long exact cohomology sequence of the pair $(M_r,U)$
twisted by $E$:
\[H^1(M_r,U)\to H^1(M_r)\to H^1(U)\stackrel{\delta}{\to}H^2(M_r,U)
\to H^2(M_r).\]

\begin{remark}
The class of a closed $1$-form $\omega$ on $U$ is contained 
in the image of $H^1(M_r)$
if and only if $\omega$ can be extended to a closed $1$-form on $M_r$.
This happens because the restriction map $\cun(M_r)\to \cun(U)$ 
is surjective.
\end{remark}

We claim that the first and last maps are zero. Indeed, a compactly-supported
form is in particular $L^2$, hence it has a smooth $L^2$ primitive by
Proposition 
\ref{cohoegze}. Thus its cohomology class on $M_r$ is zero. The long exact 
sequence therefore simplifies to
\begin{equation}\label{clos}
0\to H^1(M_r)\stackrel{i^*}{\to} H^1(U)\stackrel{\delta}{\to}H^2(M_r,U)\to 0
\end{equation}
where $i^*$ is the restriction map.
The bundle with connection $(E,D)$ does not preserve the natural
hermitian metric on $E=T_\C M$. The dual of $(E,D)$ is isomorphic to 
$(E,\overline{D})$ where $\overline{D}$ is the complex conjugate of $D$ from
Eq.\ \eqref{D}. This is isomorphic to $(E,D)$ (as real bundles)
via complex conjugation.
Thus $(E,D)$ is isomorphic to its dual.
Hence Poincar\'e duality gives 
\[H^k(M_r)\cong H^{3-k}(M_r,U)^*.\]
For $k=1$, it follows from \eqref{clos} that the (real) dimensions satisfy
\begin{equation}\label{jd}
\dim H^1(M_r)=\dim H^2(M_r,U)=\frac{\dim H^1(U)}{2}.
\end{equation}

Let us introduce the following notation for the trivalent graph $M_s$:
\begin{itemize}
\item $a$ is the number of complete geodesics (i.e., lines 
without vertices);
\item $b$ is the number of half-lines;
\item $c$ is the number of closed geodesics (loops);
\item $v$ is the number of vertices;
\item $l$ is the number of segments;
\item $g_1, \cdots, g_N$ are
the genera of the connected components of $\dr M$;
\item $n_i$ is the number of singular points on the $i$-th component of $\partial_M$.
\end{itemize}

\begin{lemma} \label{lm:dim}
The dimension of $H^1(U, E)$ equals 
\[12\sum_{i=1}^N (g_i-1)+12a+8b+4l+2k,\]
where $k$ is the 
number of independent Killing fields on $U$.
\end{lemma}

\begin{proof}
Let $S$ be the boundary of a tubular neighborhood of the singular graph $M_s$, 
viewed inside the manifold with boundary $N$ from the paragraph following 
Definition \ref{df:coco}.
The surface with boundary $S$ has $c$ connected components homeomorphic
to a torus, and $2a+b$ boundary components. 
The closed surface $\Sigma:=\partial U=\partial V$ is obtained from $\dr M$
by removing small disks around the singular points and gluing the remainder 
with $S$ 
along their common boundary circles. We use the fact that $\Sigma$ is a 
deformation retract of $U$, so they have the same (twisted and un-twisted) 
Betti numbers.
From the Mayer-Vietoris sequence,
the (untwisted) Euler characteristic of $\Sigma$ is:
\[\chi(\Sigma) = \sum_{i=1}^N \chi(\Sigma_i) +\chi(S) =\sum_{i=1}^N
  (2-2g_i-n_i) - v\]
(it is easy to see, again from the Mayer-Vietoris sequence, that $\chi(S)=-v$).
Note also the combinatorial identities
\begin{align*}
\sum n_i=2a+b,&&3v=2l+b.
\end{align*}

By lemma \ref{lkf}, the (twisted) Betti number $h^0(U)$ equals $k$, the number
of Killing vector fields on $U$. 
By Poincar\'e duality, 
since $(E,D)$ is isomorphic to $(E,\overline{D})$, we also have $h^2(\Sigma)
=h^0(\Sigma)$, therefore
$h^0(U)=h^2(U)=k$. The claim follows from the formula
\begin{equation}\label{ideul}
\chi(\Sigma,E)=\dim(E)\chi(\Sigma)
\end{equation}
(where $\chi(\Sigma,E)$ is the twisted Euler characteristic)
and from the fact that $\dim E=6$.

Eq.\ \eqref{ideul} (which is well-known)
is proved as follows: the complex
bundle $E\to\Sigma$ is flat, so its Chern character vanishes,
hence $E$ represents a torsion class in $K$-theory. This means that
$aE\oplus \C^b$ is trivial for some $a>0, b\geq 0$. Endow $\C^b$ with the
trivial 
connection, and $aE$ with the direct sum connection. By definition,
$\chi(\Sigma,aE)=a\chi(\Sigma,E)$ while $\chi(\Sigma,\C^b)=2b\chi(\Sigma)$.
Now deform the connection on $aE\oplus \C^b$ to the trivial connection. 
The Euler characteristic is constant (the index of an elliptic complex 
is always homotopy-invariant). At the end of the deformation
we get $\chi(aE\oplus \C^b)=(a\dim(E)+2b)\chi(\Sigma)$ from which 
\eqref{ideul} follows.
\end{proof}

\subsection{The local structure of the variety of representations}
\label{ssc:local}

To prove Theorem \ref{tm:existence} we need to go from an understanding 
of infinitesimal deformations, in terms of 
$H^1(M_r,E)$, to a statement on small deformations. This is based on
the inverse function theorem, applied to a natural function sending a
hyperbolic metric with cone singularities -- or more generally a 
representation of the fundamental group of $M_r$ in $PSL(2,\C)$ -- to
the induced conformal structure at infinity and cone angles. 
In this respect it is necessary to prove that the variety of 
representations of $\pi_1(M_r)$ is a smooth manifold in the neighborhood
of the holonomy representation of a convex co-compact hyperbolic metric with 
particles. 

We outline here an argument from \cite{HK,bromberg2}, also related to
earlier work of Thurston \cite{thurston-notes,culler-shalen} which can 
easily be extended to our context. This argument also provides the
dimension of the variety of representations, obtained above in a different 
way in Lemma \ref{lm:dim}.

Keeping close to the notations in \cite{HK} we call $R(M_r)$ the variety
of representations of the fundamental group of $M_r$ into $PSL(2,\C)$. 
There is a scheme associated to $R(M_r)$ by the choice of a presentation
(see \cite{weil:remarks,lubotzky-magid}), we denote it by $\cR(M_r)$. 

We first recall Theorem 5.2 of \cite{bromberg2}, which was extending a similar 
statement in \cite{HK}, itself related to a fundamental result of Thurston 
\cite{thurston-notes,culler-shalen}.

\begin{theorem}\label{tm:bromberg}
Let $M$ be a compact, connected $3$-manifold with non-empty boundary
consisting of $t$ tori and higher genus surfaces. Let $\rho\in \cR(M)$
be an irreducible representation such that, if $T$ is a torus 
component of $\dr M$, the image of $\rho(\pi_1(T))$ is neither trivial
nor $\Z_2\oplus \Z_2$. If the natural map 
$$ i:H^1(M,\dr M; E)\mapsto H^1(M; E) $$
is zero, then, in the neighborhood of $\rho$, $\cR(M)$ is a smooth
manifold of complex dimension $t-3\chi(M)+3$.
\end{theorem}

We will also need a close analog of Proposition 5.3 of Bromberg's paper 
\cite{bromberg2}, which is proved in the same way. 

\begin{prop} \label{pr:bromberg}
Let $M$ be a convex co-compact hyperbolic 3-manifold and let $\rho$ be
its holonomy representation.
\begin{enumerate}
\item The restriction of $\rho$ to each end is irreducible.
\item $\rho$ is irreducible. 
\item Let $T$ be a tubular neighborhood of a closed curve which is a
connected component of the singular locus $M_s$. Then $\rho(\pi_1(T))$ is
infinite and non-parabolic. 
\end{enumerate}
\end{prop}

\begin{proof}[Sketch of the proof]
The first two points were already noted in subsection \ref{ssc:variety}. 
The last point follows from the fact that the
holonomy of any closed curve in $T$ which is not a multiple of a 
meridian has to be loxodromic.
\end{proof}

We can now check that Theorem \ref{tm:bromberg} can be applied in
our context. The tori boundary components of $N$ correspond to the 
closed curves in the singular locus of $M$, so that the condition that
$\rho(\pi_1(T))\not\in \{ 1, \Z_2\oplus Z_2\}$ follows from point (3) of
the previous proposition. 
Moreover, Theorem \ref{tm:rigidity} immediately implies that the map $i$
appearing in the statement of Theorem \ref{tm:bromberg} is zero: 
if $\omega\in H^1(N,\dr N; E)$ then $\omega$ does not changer either
the conformal structure at infinity or the angles at the cone
singularities, so that $\omega=0$. 
So we can conclude that, if $M$ is a 
convex co-compact hyperbolic manifold with particles, then
$R(\pi_1(M_r))$ is a smooth manifold in the neighborhood
of the holonomy representation $\rho$ of $M_r$. 

Note that the dimension of $R(\pi_1(M_r))$ from Theorem \ref{tm:bromberg}
is not the same as the dimension of the space of deformations of $M$
among hyperbolic cone-manifolds (plus the dimension of $PSL(2,\C)$, 
because some representations close
to $\rho$ in $R(M_r)$ do not come from hyperbolic metrics with 
cone singularities -- this happens precisely when the holonomy of
the meridian of the tubular neighborhood of a cone singularity
is not elliptic but has a translation component (this condition appears
in the definition of $\cR_{cone}(M_r)$ just before Theorem 
\ref{tm:existence}).

\subsection{Proof of Theorem \ref{tm:existence}}

We now have the tools necessary to prove Theorem \ref{tm:existence}.
First note the following

\begin{lemma}\label{kc}
There do not exist non-zero Killing vector fields on the connected components
of $U$ others than those corresponding to the closed geodesic loops, 
therefore the number $k$ in Lemma \ref{lm:dim} is $2c$. 
\end{lemma}

\begin{proof}
Let $\kappa$ be such a Killing field on a connected component of $U$ which
contains a geodesic line or half-line from $M_s$. Then $\kappa$
would  have an extension as a holomorphic vector field $v_\kappa$
on the boundary at infinity of $M$. Moreover, since the angles at the
singular arcs are less than $\pi$, any Killing field has to behave, near each
singular arc, as a Killing field in $V_\alpha$ with axis $\Delta_0$ --- indeed
the only Killing fields on $V_\alpha, 0<\alpha<\pi$, are induced by Killing
fields on $H^3$ with axis $\Delta_0$. It follows that $v_\kappa$ has zeros at
the singular points of $\dr_\infty M$, i.e., at the endpoints of the singular
arcs. 

Consider a connected component $\dr_{\infty,0}M$ of $\dr_\infty M$, and the
corresponding connected component $\dr_0 CC(M)$ of the boundary of the convex
core of $M$. $\dr_0 CC(M)$ is ruled and convex, therefore hyperbolic, outside
its intersections with the singular locus of $M$, where it has singular points
of singular curvature less than $2\pi$. It follows from the Gauss-Bonnet
theorem that if $\dr_0 CC(M)$ is a torus, it intersects at least one singular
arc, while if it is a sphere, it intersects at least 3 singular arcs.

The vector field $v_\kappa$ considered above is holomorphic, and it has at
least 3 zeros on $\dr_{\infty, 0}M$ if $\dr_{\infty, 0}M$ is a sphere, and
at least one if $\dr_{\infty, 0}M$ is a torus. Therefore it vanishes. 
So $v_\kappa$ vanishes on $\dr_\infty M$, and it follows that $\kappa$ is
zero. 

For a component which contains a trivalent vertex, it is geometrically obvious
that there are no local isometries near that vertex.
\end{proof}

Let $\cD$ be the space of data appearing in Theorem \ref{tm:existence}, with
in addition, for each singular arc (either a segment, a circle,
a half-line or a line), 
a number corresponding to the translation component of the holonomy 
along that singularity. Thus:
\[\cD = \left(\Pi_{i=1}^N \cT_{g_i, n_i}\right) \times (\R_+)^{a+b+c+l} 
\times \R^{a+b+c+l}~,\]
where, for each $i\in \{ 1,\cdots, N\}$, $n_i$ is the number of endpoints of
the singular arcs on $\dr_iM$, the factor $(\R_+)^{a+b+c+l}$ corresponds to the
angles around the singular arcs, and the term $\R^{a+b+c+l}$ corresponds to the
translation component along the singular arcs of the corresponding holonomy. 
The factors $\cT_{g_i, n_i}$ contain the conformal
structure at infinity on $\dr_iM$, with marked points corresponding to the
endpoints of the singular arcs.

\begin{remark}
$\dim(\cD)=6\sum_{i=1}^N (g_i-1)+6a+4b+2c+2l$.  
\end{remark}

\begin{proof}
For each $i\in \{ 1,\cdots, N\}$, $\dim(\cT_{g_i, n_i})=6g_i-6+2n_i$, so the
formula follows from the fact that $\sum_{i=1}^N n_i=2a+b$ because each singular
line has two endpoints, each singular half-line has one endpoint on $\partial M$, 
while the edges and the circles have no such endpoint.
\end{proof}

From this remark, from Lemma \ref{lm:dim}, Lemma \ref{kc} and Eq.\ \eqref{jd}, we deduce
the following

\begin{cor}
The dimensions of $H^1(M_r, E)$ and of $\cD$ are equal. 
\end{cor}

Let $g$ be a convex co-compact hyperbolic singular metric on $M$, as in
Theorem \ref{tm:existence}, and let $c(g)$ be the induced element of $\cD$. By
definition the last term in $c(g)$, in $\R^a$, is equal to $0$, since 
$g$ is a cone-manifold (so that the translation component of the holonomy
is $0$ for each singular arc). Each element of $H^1(M_r)$ induces an infinitesimal
variation of the holonomy of $(M_r,g)$, and therefore an element of the tangent
space $T_{c(g)}\cD$ of $\cD$ at $c(g)$, and this defines a linear map:
$$ \gamma:H^1(U,E) \rightarrow T_{c(g)}\cD~. $$

Since $H^1(M_r, E)$ and $\cD$ have the same dimension, and 
$\gamma\circ i^*:H^1(M_r,E)\to T\cD$ is injective by Theorem
\ref{tm:rigidity}, it follows that it is surjective. 
But $H^1(M_r,E)$ is the tangent space of $R(M_r)/PSL(2,\C)$ at the holonomy
representation of the hyperbolic metric with cone singularities 
considered, and $\gamma$ is the differential
of the (smooth) map sending a holonomy representation of a hyperbolic metric with cone 
singularities to its cone angles and conformal structure at infinity. It is
therefore possible to apply the inverse function theorem, which yields
Theorem \ref{tm:existence}. 

\appendix

\section{Convex subsets in hyperbolic manifolds with particles}

\subsubsection*{Definitions, outline.}

This appendix contains some basic information on the geometry of
convex subsets in hyperbolic manifolds with particles. The term ``convex''
should be understood here as in Definition \ref{df:convex}: a non-empty
subset $K$ is convex if any geodesic segment with endpoints in $K$ is
entirely contained in $K$. 

It follows directly from this definition that the intersection of two
convex subsets is either empty or convex. Our main goal here is to show that,
under weak topological assumptions on $M$,
the intersection of two non-empty convex subsets cannot be empty. It will
follow that it is possible to define the convex core of a hyperbolic manifold
with particles, and we will then point out some of its elementary properties. 

As in the body of the paper we consider here a hyperbolic manifold with
particles $M$, and denote by $M_r$ and by $M_s$ its regular and singular set,
respectively. By definition $M_s$ is a finite graph, and the angle at each of
its edges is less than $\pi$.

\subsubsection*{Links of points in $M$.}

It is useful to consider the set of unit vectors based at a point of
$M$. For regular points of $M$ this is just the unit tangent bundle, however
for singular points this notion is more interesting.

\begin{defi}
Let $x\in M$, we call $L_x(M)$ (or simply $L_x$) the set of unit vectors at
$x$. $L_x(M)$ is the {\it link} of $M$ at $x$.
\end{defi}

Note that a unit vector can be defined (in a general setting) as the speed
at $x$ of a geodesic ray starting at $x$ with velocity $1$; 
two unit vectors are identical if
the corresponding geodesic rays are equal in an interval containing $0$.
There is a natural distance on $L_x$, defined by the angle between two unit
vectors. It follows from the definition of a cone-manifold that $L_x$, with 
this distance, is a spherical surface with cone singularities. The cone points
correspond to the singular segments containing $x$, and the angle at those
cone points in $L_x$ is equal to the angle at the corresponding singular
arc in $M$.

There is a particular kind of spherical cone-manifolds which plays an 
important role here. Let $\theta\in (0,\pi)$, consider the universal cover of
the complement of the two ``poles'' in the unit sphere, and then its quotient
by the rotation of angle $\theta$ fixing the two poles. This quotient is
denoted here by $S^2_\theta$, it is a spherical surface with two cone
singularities where the angle is equal to $\theta$.

\begin{prop}\label{locd}
Let $x\in M$, then:
\begin{itemize}
\item if $x\in M_r$ then $L_x$ is isometric to the unit sphere $S^2$,
\item if $x\in M_s$ is contained in the interior of a singular edge $e$, 
then $L_x$ is isometric to $S^2_\theta$, where $\theta$ is the angle at $e$,
\item if $x\in M_s$ is a vertex, then $L_x$  is isometric to the metric space
  obtained by gluing along their common boundary two copies
  of a spherical triangle with acute angles.
\end{itemize}
\end{prop}

\begin{proof}
If $x\in M_r$ the statement is quite obvious. In the second case the result
can be obtained from the definition of the hyperbolic metric in the 
neighborhood of a point in the interior of a singular segment. 

In the third
case the link of $x$ is by definition a spherical metric with cone
singularities. Moreover, the angle at each cone point is less than $\pi$, so
that the corresponding singular curvature is larger than $\pi$. So it follows
from the Gauss-Bonnet theorem that there are at most 3 cone points. But it
follows from a theorem of Alexandrov (see \cite{alex,luo-tian}) that such a
metric is the double cover of a spherical triangle -- this 
can also be proved
directly, without reference to Alexandrov's much more general theorem. 
Finally, since the angle
at each cone point is less than $\pi$, the spherical triangle has acute
angles. 
\end{proof}

Note that in the third case the angles at the three singular segments
arriving at $x$ are twice the angles of a spherical triangle (with acute
angles). These angles are equal to the edge lengths of the dual spherical
triangle (for the polar duality in the sphere) so they satisfy the triangle
inequality. It follows that the angles at the three singular arcs containing
$x$ also satisfy the triangle inequality, and the same line of reasoning shows
that any triple of angles in $(0,\pi)$ satisfying the triangle inequality can
be realized in this manner.

\subsubsection*{The link of a convex subset at a point.}

Now let $K$ be a convex subset of $M$. 

\begin{defi}
Let $x\in K$, we call $L_x(K)$ the set of unit vectors $v\in L_x(M)$ such that
the geodesic ray starting from $x$ in the direction of $v$ is contained (on
some interval containing $0$) in $K$. $L_x(K)$ is the link of $K$ at $x$.
\end{defi}

Clearly $L_x(K)=\emptyset$ when $x$ is not contained in $K$, while
$L_x(K)=L_x$ when $x$ is contained in the interior of $K$. The most
interesting case is when $x\in \dr K$, then $L_x(K)$ is a subset of $L_x(M)$. 
This subset is almost always geodesically convex in the following sense.

\begin{lemma} \label{lm:pi}
Let $x\in \dr K$, and let $\gamma$ be a geodesic segment in $L_x(M)$ of length
less than $\pi$, with endpoints in $L_x(K)$. Then
$\gamma$ is contained in $L_x(K)$. 
\end{lemma}

\begin{proof}
Note that no cone point of $L_x(M)$ is contained in $\gamma$, except perhaps
at its endpoints. 
Let $\epsilon>0$ be small enough, and let $\Omega_\epsilon$ be the union 
of the geodesic segments of length $\epsilon$ starting from $x$ in the
directions of $\gamma$. Then $\Omega_\epsilon$ is a plane sector of angle less
than $\pi$ at $x$.

Let $s\in (0,\epsilon)$, consider the geodesic segment $c_s$ in
$\Omega_\epsilon$ with endpoints the points at distance $s$ from $x$ in
$\Omega_\epsilon$ in the segments starting from $x$ in the direction of the
endpoints of $\gamma$.
By definition of a convex subset, $c_s$ is contained in
$K$. This shows that a neighborhood of $x$ in $\Omega_\epsilon$ is contained
in $K$, and therefore that $\gamma$ is contained in $L_x(K)$.
\end{proof}

\begin{cor}
If $L_x(K)$ has dimension $2$, then it has locally convex boundary in $L_x(M)$.
\end{cor}

\begin{cor} \label{cr:point}
Let $K\subset M$ be a convex subset, which is not reduced to one point. 
Then $\dr K$ contains no vertex of $M$.   
\end{cor}

\begin{proof}
Suppose that $v$ is a vertex of $M$, then $L_v(M)$ is obtained by ``doubling''
a spherical triangle with acute angles, we call $c_1, c_2$ and $c_3$
its cone singularities. Suppose now that $v\in \dr K$, and that $K$ is not reduced to
one point, so that $L_v(K)\neq \emptyset$. Since $L_v(K)$ is the double
of a spherical triangle with acute angles, its diameter is less than
$\pi$, it then follows from Lemma \ref{lm:pi} that $L_v(K)$ is connected.

Note that $L_v(K)$ cannot be reduced to only one point which is a 
cone singularity of $L_v(M)$. Indeed, suppose for instance that 
$L_v(K)=\{ c_1\}$, let $\gamma$ be a simple loop based at $c_1$ with 
$c_2$ on one side and $c_3$ on the other. Since the cone angles are less
than $\pi$, $\gamma$ can be deformed (in the complement of the singular
points) to a geodesic loop based at $c_1$.
A standard argument in the geometry of spherical surfaces shows that 
this geodesic loop has length less than $\pi$ (this uses the
fact that the cone angles are less than $\pi$). So $\gamma\subset L_v(K)$
by Lemma \ref{lm:pi}, and therefore $L_v(K)$ contains points of $L_v(M)$
other than $c_1$.

Let $x$ be a point of $L_v(K)$ which is not a cone
singularity of $L_v(M)$. Let $\gamma_1, \gamma_2$ and $\gamma_3$ be simple loops
based at $x$, and going around $c_1, c_2$ and $c_3$,
respectively. Since the angles at the $c_i$ are less than $\pi$, the
curves $\gamma_i$ can be deformed (in the complement of the singular
points) to minimizing geodesic loops (based
at $x$), and those curves, being minimizing, are disjoint. 

As already used above, the lengths of the $\gamma_i$ are less than $\pi$,
so Lemma \ref{lm:pi}
shows that the $\gamma_i$ are contained in $L_v(K)$. A simple convexity
argument then shows that $L_v(K)$ contains a neighborhood of $x$, 
so that $L_v(K)$ is non-degenerate (it has dimension 2). 

The complement of $\gamma_1\cup\gamma_2\cup \gamma_3$ is composed of
four topological disks, three containing one of the cone singularities
of $L_v(M)$, and the last one not containing any. Since 
$L_v(K)$ has locally convex boundary, the Gauss-Bonnet formula shows
that it has positive Euler characteristic, so $L_v(K)$
has to contain at least 3 of the 4 disks in the complement of 
$\gamma_1\cup\gamma_2\cup \gamma_3$ because it has genus $0$. 
But then $L_v(K)$ contains at
least two cone singularities, and the sum of their singular curvatures 
is more than $2\pi$, so that the Euler characteristic of $L_v(K)$ has
to be at least $2$ by the Gauss-Bonnet formula. So $L_v(K)=L_v(M)$,
and this contradicts the fact that $v\in \dr K$.
\end{proof}

The same kind of arguments can be used to understand the link of $K$
at a boundary point which is contained in a singular segment, but is not
a vertex of the singular set of $M$.

\begin{cor} \label{cr:a7}
Let $v\in \dr K$ be contained in a singular arc of $M$ with angle $\theta$, 
then $L_v(K)$ is a subset of $L_v=S^2_\theta$.  If $L_v(K)$ contains a point
which is not one of the cone singularities of $S^2_\theta$, then 
\begin{itemize}
\item either $L_v(K)$ is non-degenerate (i.e., it has dimension 2) and it 
contains exactly one of the cone singularities of $S^2_\theta$,
\item or $L_v(K)$ is a closed geodesic, and it contains no cone singularity.
\end{itemize}
\end{cor}

\begin{proof}
Let $x\in L_v(K)$ be a point which is not one of the cone points of $L_v(M)$. 
Let $\gamma$ be a simple loop based at $x$, not homotopically trivial
in the complement of the cone points of $L_v(M)$. Then $\gamma$ can be deformed
to a geodesic loop, of length less than $\pi$ (because the cone angles of
$L_v(M)$ are less than $\pi$) so it is contained in $L_v(K)$.

This geodesic loop can be a closed geodesic, in this case it can be equal to
$L_v(K)$. This corresponds to the second case in the statement of the corollary.
We now suppose that we are not in this case. Then 
$L_v(K)$ contains a neighborhood of $x$, so it is non-degenerate. Since $L_v(K)$
has locally geodesic boundary, it has positive Euler characteristic by the
Gauss-Bonnet Theorem, so $L_v(K)$ must contain one of the disks bounded by
$\gamma$, so one of the cone singularities of $L_v(M)$. 

But $L_v(K)$ cannot contain both cone singularities of $L_v(M)$, otherwise its
Euler characteristic would be at least $2$, again by the Gauss-Bonnet Theorem
because the sum of the singular curvatures of the cone singularities is larger
than $2\pi$. Thus $L_v(K)$ would be equal to $L_v(M)$, this is impossible
since $v\in \dr K$.
\end{proof}

\subsubsection*{The normal unit bundle.}

Here we consider a convex subset $K$ in $M$ (the definitions given here
make sense for other subsets).

\begin{defi} \label{df:n1}
Let $x\in \dr K$, the unit normal subset of $K$ at $x$, called $N^1_xK$,
is the set of points
$v\in L_x(M)$ which are at distance at least $\pi/2$ from $L_x(K)$.
\end{defi}

For instance:
\begin{itemize}
\item If $x\in M_r$, then $L_x(K)$ is a subset of $L_x$ with locally
convex boundary, and $N^1_x(K)$ is the dual of $L_x(K)$.
\item In particular, if
$\dr K$ is smooth at $x$, then $N^1_xK$ has only one point, which is the
unit normal of $\dr K$ at $x$.
\item If $K=\{ x\}$, then $L_x(K)=\emptyset$ and $N^1_xK=L_x$.
\item If $x\in K\subset M_s$ and $K$ contains a segment of $M_s$ around
$x$, then  $L_x(K)$ is made of the two singular points of $L_x$, and 
$N^1_xK$ is the ``equator'' of $L_x$ (the set of points at distance 
$\pi/2$ from both singular points).  
\end{itemize}

\begin{remark} \label{rk:closest}
Let $y\in M\setminus K$, and let $x$ be a point in $\dr K$ such that
$d(x,y)=d(K,y)$. Let $c$ be a minimizing geodesic segment between $x$ and
$y$. Then the unit vector at $x$ in the direction of $c$ is contained in
$N^1_x(K)$.   
\end{remark}

\begin{proof}
Clearly, otherwise it would be possible to find another point of $K$, close to
$x$, closer to $y$ than $x$.  
\end{proof}

We call $N^1(K)$ the disjoint union of the sets $N^1_x(K)$ over the points
$x\in \dr K$.

\subsubsection*{The normal exponential map.}

Let $x\in M$, let $v\in L_x$ and let $t\in \R_+$, we denote by $\exp_x(tv)$
the point of $M$ which is at distance $t$ from $x$ on the geodesic ray
starting from $x$ with speed $v$. Note that, given $x$ and $v$, $\exp_x(tv)$
is well-defined for $t$ small enough (if $v$ is a regular point of $L_x$, 
until the geodesic ray starting
from $x$ in the direction of $v$ arrives at the singular set of $M$, and, if
$v$ is a cone point of $L_x$, until that geodesic arrives at a vertex of
$M$). 

Note that it is not clear at this point that $\exp$ is defined at all points
of $N^1(K)\times \R_+$ since some geodesic rays could run into the singular
set of $M$. We will see below that this can not happen. In the meantime we
call $R$ the length of the smallest geodesic segment where this phenomenon
happens; thus it will be shown below that $R=\infty$.

\begin{lemma} 
The map: $\exp:N^1K\times (0,R)\rightarrow M$ is a
homeomorphism from $N^1K\times (0, R)$
to the set of points at distance less than $R$ from $K$ in $M\setminus K$.
It sends the complement of the points $(x,v,t)$ where $x\in
\dr K$ and $v$ is a singular point of $L_x$ to the complement of the singular
locus in $M\setminus K$.
\end{lemma}

\begin{proof}
By construction the restriction of $\exp$ to $N^1(K)\times (0,r)$ is a
homeomorphism onto its image for $r$ small enough. Moreover Remark
\ref{rk:closest} shows that its image is exactly the set of points at
distance less than $r$ from $K$. The shape operators of the surfaces
$\exp(N^1(K)\times \{ s\})$, for $s\in (0,r)$, satisfy a Riccati equation,
and an argument which is classical in hyperbolic geometry shows that these
surfaces are locally convex. It also follows from the definition that they are
orthogonal to the singular locus.

Suppose that $\exp$ is not injective on $N^1(K)\times (0,R)$. 
Let $r_M$ be the supremum of the $r\in (0,R)$ such that the restriction of
$\exp$ to $N^1(K)\times (0,r)$ is injective. Since $\exp$ remains a local
homeomorphism at $r$, there are two points
$(x,v), (x',v')\in N^1(K)$ such that $\lim_{r\rightarrow r_M}\exp_x(rv) =
\lim_{r\rightarrow r_M}\exp_{x'}(rv')$. But then the set:
$$ \exp_x([0,r_M]v)\cup \exp_{x'}([0,r_M]v') $$
is a geodesic segment (otherwise an intersection would appear before $r=r_M$)
with endpoints in $K$ but which is not contained in $K$,
a contradiction. 
\end{proof}

Now suppose that $R<\infty$, then there is a geodesic segment $c$ of length
$R$ starting from a point $x\in \dr K$, with direction given by a vector $v\in
N^1_x(K)$, and ending at a point $y\in M_s$. Moreover, either $c$ is contained
in $M_r$ (except for its endpoints) and $y$ is in an arc $e$ of $M_s$, 
or $c$  is contained in an
arc of $M_s$, and $y$ is a vertex of $M$. Let $w$ be the unit vector at $y$ in
the direction of $c$. 

In the first case, $w$ is a point of $L_y$ at
distance $\pi/2$ of both cone points of $L_y$. Let $\gamma$ be a geodesic
ray starting from $y$ in a direction $w'$ which is not one of the cone points
of $L_y$. Then the distance between $w$ and $w'$ in $L_y$ is less than 
$\pi/2$, so the derivative of the distance to $K$ is negative along $\gamma$. 
So for $r\in (0,R)$
close enough to $R$, the set of points at distance less than $r$ from $K$
contains the complement of the $M_s$ in a neighborhood of $y$, so it  
does not retract on $K$. This contradicts the previous lemma.

In the second case, $y$ is a vertex of $M_s$, and all points of $L_y$ are at
distance less than $\pi/2$ from $w$ (this follows from the description of
$L_y$ as obtained by gluing two copies of a spherical triangle with acute
angles). So $y$ is a local maximum of the distance to $K$, and this yields
again a contradiction with the previous lemma. So $R=\infty$. This argument
shows the following statement.

\begin{lemma} \label{lm:ends}
The map $\exp:N^1(K)\times (0,\infty)\rightarrow M\setminus K$ is a
homeomorphism.
It sends the complement of the points $(x,v,t)$ where 
$x\in \dr K$ and $v$ is a singular point of $L_x$ to the regular set
of $M\setminus K$. 
\end{lemma}

\subsubsection*{A global description of convex subsets.}

It follows from the previous lemma that no ``accident'' occurs in the
map $\exp:N^1(K)\times (0,R)\rightarrow M$. Therefore:
\begin{itemize}
\item all vertices of $M_s$ are contained in $K$,
\item for each point $y\in M_s$ outside $K$, there is only one geodesic segment
  minimizing the distance from $y$ to $K$, and it is contained in an edge of
  $M_s$. 
\end{itemize}
The map $\exp$ defines a homeomorphism 
from $\Sigma\times (0,\infty)$,
where $\Sigma$ is a (non connected) closed surface, to $M\setminus K$.
Here $\Sigma$ corresponds to $N^1K$ so that it is only a $C^0$ surface.
Moreover $M_s\setminus K$ is the image of $S\times (0,\infty)$, where 
$S$ is a finite subset of $\Sigma$.

The following is a rather direct consequence. 

\begin{lemma} \label{lm:closed-geod}
Let $K\subset M$ be a non-empty convex subset. Then $K$ contains
\begin{enumerate}
\item all vertices of the singular set $M_s$,
\item all closed geodesics in $M$.
\end{enumerate}
\end{lemma}

\begin{proof}
The first point is a direct consequence of Lemma \ref{lm:ends}. 
For the second point we use a ``trick'' already used before, and define
$u:M\rightarrow \R$ by 
$$ u(x) = \sinh(d(x,K))~. $$
It is then known that $u$ satisfies on $M\setminus K$ the inequality
$$ \mbox{Hess}(u) \geq u g~, $$
where $g$ is the hyperbolic metric on $M$. The reason for this inequality
is that, in $H^3$, the $\sinh$ of the distance to a totally geodesic plane
satisfies the equality $\mbox{Hess}(u)=ug$, the inequality then follows 
from a simple local argument using Lemma \ref{lm:ends} and the local 
convexity of the boundary of $K$. Note that the inequality should be
understood in a distribution sense if the boundary of $K$ is not smooth.

If $\gamma$ is a closed geodesic in $M$, parametrized at velocity $1$, 
then it follows that $(u\circ \gamma)''\geq u\circ \gamma$ (in a distribution
sense), which is
clearly impossible by the maximum principle unless $u\circ \gamma=0$. 
This shows the second point in the lemma.
\end{proof}

\subsubsection*{Convex subsets have non-empty intersection.}

We are now ready to obtain the result announced at the beginning of this
appendix. 

\begin{lemma} \label{lm:intersection}
Suppose that either $M_s$ has a vertex or $\pi_1(M)\neq 0$. 
Let $K, K'$ be two non-empty compact convex subsets of $M$, then 
$K\cap K'\neq\emptyset$.
\end{lemma}

\begin{proof}

Suppose first that $M_s$ has at least one vertex $v$, then Lemma 
\ref{lm:closed-geod} shows that both $K$ and $K'$ contain $v$, so
$v\in K\cap K'$. 

Suppose now that $\pi_1(M)\neq 0$, let $\gamma$ be a closed curve in a
non-trivial element of $\pi_1(M)$. Let $(\gamma_n)_{n\in \N}$ be a 
minimizing sequence in the homotopy class of $\gamma$, that is, a sequence
of curves homotopic to $\gamma$ such that the length of $\gamma_n$
converges, as $n\rightarrow \infty$, to the infimum of the lengths
of curves homotopic to $\gamma$. It follows from Lemma \ref{lm:ends}
-- and from the form of the ends of $M$ -- that the $\gamma_n$
remain at bounded distance from $K$. So, after extracting a subsequence,
$\gamma_n$ converges to a closed geodesic $\gamma_\infty$. Lemma 
\ref{lm:closed-geod} shows that $\gamma_\infty$ is contained 
in both $K$ and $K'$, and the result follows.
\end{proof}

It would be useful to weaken the hypothesis of this lemma by supposing 
only that $\pi_1(M_r)\neq 0$. Such an extension might be true, but some care
is required. It is quite possible that, if $\pi_1(M_r)\neq 0$ (and $M$
contains a non-empty compact convex subset) then $M_r$ contains
a closed geodesic. It is however {\it not} true that any non-trivial element of 
$\pi_1(M_r)$ can be realized as a closed geodesic. 

If $M=V_\alpha$, for some $\alpha\in (0,\pi)$, then Lemma \ref{lm:intersection}
does not apply: if for instance $K$ and $K'$ are each reduced to one point 
in the singular locus, then $K\cap K'$ could be empty. There are also of course
hyperbolic manifolds (with or without singularities) which do not contain any
non-empty compact convex subset. 

A more interesting example is obtained from a hyperbolic metric $h$ with four 
cone singularities of angle less than $\pi$ on the sphere $S^2$. One can then
consider the warped product metric 
$$ dt^2 + \cosh^2(t) h $$
on $S^2\times \R$, it is easily seen to be a complete hyperbolic metric with
four cone singularities along infinite lines. It contains $S^2\times \{ 0\}$
as a compact convex subset. Theorem \ref{tm:existence} shows that this
example can be deformed, by changing the conformal structure at infinity. 
However, $\pi_1(S^2\times \R)=0$, so that Lemma \ref{lm:intersection} does not 
apply. It would be desirable to have a more general statement including this
example. 

\subsubsection*{The boundary of the convex core.}

Lemma \ref{lm:intersection} shows that it is possible to define the convex
core $CC(M)$ of $M$ as the smallest non-empty convex subset in $M$. By the
considerations above, $CC(M)$ contains all the vertices of $M_s$, and
$M\setminus CC(M)$ is the disjoint union of ``ends'', each of which is
homeomorphic to the product of a closed surface by an interval. 

\begin{lemma} \label{lm:cc-ortho}
The boundary of $CC(M)$ is a surface orthogonal to the singular locus.
\end{lemma}

\begin{proof}
Let $x\in M_s\cap \dr CC(M)$. Then $x$ is contained in a singular arc $e$ of
$M_s$, let $\theta$ be the angle around $e$.
By construction, $L_x(CC(M))$ is a subset of 
$L_x=S^2_\theta$. We have seen in Corollary \ref{cr:a7} that 
\begin{itemize}
\item either $L_x(CC(M))$ is a closed geodesic, 
\item or $\dr L_x(CC(M))\subset S^2_\theta$ is a locally 
convex curve and $L_x$ contains the ``south pole'' $p_S$ of
$S^2_\theta$, i.e., the image in $S^2_\theta$ of the ``south pole'' in $S^2$. 
\end{itemize}
In the first case $CC(M)$ is a totally geodesic surface in the neighborhood
of $x$, and a simple connectedness argument shows that it is everywhere
totally geodesic, and thus orthogonal to the singular locus. We therefore
consider the second case.

Consider the function $\phi$ defined on $S^2_\theta$ as the distance to
$p_S$. $\phi(p_S)=0$, while $\phi(p_N) = \pi$, where $p_N$ is the
``north pole'' of $S^2_\theta$. Let $y\in \dr L_x(CC(M))$ 
be the point where $\phi$
attains its minimum. Consider the geodesic segment $\gamma:(-l,l)\rightarrow
S^2_\theta$, parametrized at speed $1$, such that $(\phi\circ \gamma)'(0)=0$,
where $l$ is chosen to be maximum under the condition that
$\gamma$ is embedded. Then $\lim_{-l}\gamma = \lim_l\gamma$, so that the 
closure of $\gamma((-l,l))$ is a closed curve, which is geodesic except at one
point. Note that $\gamma((-l,l))$ is simply the projection to $S^2_\theta$ of
a geodesic segment of length $2l$ in $S^2$. 
The local convexity of $\dr L_x$ shows that it remains ``under''
$\gamma((-l,l))$. It follows that: 
\begin{itemize}
\item If $\phi(y)>\pi/2$, then the restriction of $\phi$ to $\gamma((-l,l))$
attains a strict maximum at $y$, so the restriction of $\phi$ to $\dr L_y$ 
is also maximal at $y$, and takes strictly lower values at other points. 
This clearly contradicts the definition of $y$, so this case can be
eliminated. 
\item If $\phi(y)<\pi/2$, then the fact that $\theta<\pi$ implies that the
restriction of $\phi$ to $\gamma((-l,l))$ is negative, so the restriction of
$\phi$ to $\dr L_x(CC(M))$ is also negative. This contradicts the definition of
$CC(M)$ as the smallest convex subset in $M$, because if would then be
possible to reduce $CC(M)$ by ``cutting out'' the part above a plane
orthogonal to the singular locus but slightly ``below'' $x$, and still get a
convex subset of $M$.
\item If $\phi(y)=\pi/2$, then $\phi$ is identically $\pi/2$ on $\gamma((-l,l))$, so
  that the restriction of $\phi$ to $\dr L_x(CC(M))$ 
  is at most $\pi/2$. The definition
  of $y$ as the point where $\phi$ is minimum thus entails that $\phi$ is
  identically $\pi/2$ on $\dr L_x(CC(M))$, 
  and this means precisely that $\dr CC(M)$ is
  orthogonal to the singular locus at $x$.
\end{itemize}
This argument shows that $\dr CC(M)$ is orthogonal to the singular locus at
$x$, as claimed.
\end{proof}

This has interesting consequences, which can be summed up as follows.

\begin{lemma} \label{lm:convex-core}
The boundary of $CC(M)$ is a ``pleated surface''. Its induced metric is
hyperbolic, with cone singularities at the intersection with the singular arcs
in $M_s$, and the angle at each such cone point is equal to the angle at the
corresponding singular arc of $M_s$. The surface $\dr CC(M)$ is ``bent'' along
a measured lamination whose support is disjoint from the cone points. 
\end{lemma}

\begin{proof}[Sketch of the proof.]
We do not give complete details of the proof, which is similar to the
corresponding situation with no ``particle'', as in \cite{thurston-notes}. 
The fact that $\dr CC(M)$ is a pleated surface away from $M_s$ is a
consequence of the fact that it is the boundary of a convex subset
of $M$ without extremal point, as in the non-singular case. 
Since $\dr CC(M)$ is orthogonal to the
singular arcs, its induced metrics has, at those intersections, cone points
with angle equal to the angle at those singular arcs in $M$. The same fact
also entails that the support of the bending lamination does not contain the
singular points (and therefore, since the angles at the cone points are less
than $\pi$, the distance from the support of the bending lamination to the
cone points is bounded away from $0$).
\end{proof}

\bibliographystyle{amsplain}

\begin{thebibliography}{10}

\bibitem{ahlfors}
L.~V. Ahlfors, \emph{Lectures on quasiconformal mappings}, D. Van Nostrand Co.,
  Inc., Toronto, Ont.-New York-London, 1966, Manuscript prepared with the
  assistance of Clifford J. Earle, Jr. Van Nostrand Mathematical Studies, No.
  10.

\bibitem{ahlfors-bers}
Lars Ahlfors and Lipman Bers, \emph{Riemann's mapping theorem for variable
  metrics}, Ann. of Math. (2) \textbf{72} (1960), 385--404.

\bibitem{alex}
A.~D. Alexandrov, \emph{Convex polyhedra}, Springer Monographs in Mathematics,
  Springer-Verlag, Berlin, 2005, Translated from the 1950 Russian edition by N.
  S. Dairbekov, S. S. Kutateladze and A. B. Sossinsky, With comments and
  bibliography by V. A. Zalgaller and appendices by L. A. Shor and Yu. A.
  Volkov.

\bibitem{bonahon-otal}
Francis Bonahon and Jean-Pierre Otal, \emph{Laminations mesur\'ees de plissage
  des vari\'et\'es hyperboliques de dimension 3}, Ann. Math. \textbf{160}
  (2004), 1013--1055.

\bibitem{cone}
Francesco Bonsante and Jean-Marc Schlenker, \emph{{AdS} manifolds with
  particles and earthquakes on singular surfaces}, {math.GT/0609116. Geom.
  Funct. Anal., to appear.}, 2006.

\bibitem{bromberg1}
K.~Bromberg, \emph{Hyperbolic cone-manifolds, short geodesics, and {S}chwarzian
  derivatives}, J. Amer. Math. Soc. \textbf{17} (2004), no.~4, 783--826
  (electronic). \MR{MR2083468}

\bibitem{bromberg2}
\bysame, \emph{Rigidity of geometrically finite hyperbolic cone-manifolds},
  Geom. Dedicata \textbf{105} (2004), 143--170. \MR{MR2057249}

\bibitem{Calabi}
E.~Calabi, \emph{On compact {Riemannian} manifolds with constant curvature,
  {I}}, AMS proceedings of Symposia in Pure Math \textbf{3} (1961), 155--180.

\bibitem{culler-shalen}
Marc Culler and Peter~B. Shalen, \emph{Varieties of group representations and
  splittings of {$3$}-manifolds}, Ann. of Math. (2) \textbf{117} (1983), no.~1,
  109--146. \MR{MR683804 (84k:57005)}

\bibitem{c-epstein}
Charles~L. Epstein, \emph{The hyperbolic {G}auss map and quasiconformal
  reflections}, J. Reine Angew. Math. \textbf{372} (1986), 96--135.
  \MR{MR863521 (88b:30029)}

\bibitem{fefferman-graham}
Charles Fefferman and C.~Robin Graham, \emph{Conformal invariants},
  Ast\'erisque (1985), no.~Numero Hors Serie, 95--116, The mathematical
  heritage of \'Elie Cartan (Lyon, 1984).

\bibitem{GHL}
Sylvestre Gallot, Dominique Hulin, and Jacques Lafontaine, \emph{Riemannian
  geometry}, third ed., Universitext, Springer-Verlag, Berlin, 2004.
  \MR{MR2088027}

\bibitem{graham-lee}
C.~R. Graham and J.~M. Lee, \emph{Einstein metrics with prescribed conformal
  infinity on the ball}, Adv. Math. \textbf{87} (1991), 186--225.

\bibitem{HK}
Craig~D. Hodgson and Steven~P. Kerckhoff, \emph{Rigidity of hyperbolic
  cone-manifolds and hyperbolic {Dehn} surgery}, J. Differential Geom.
  \textbf{48} (1998), 1--60.

\bibitem{kostant}
Bertram Kostant, \emph{Holonomy and the {L}ie algebra of infinitesimal motions
  of a {R}iemannian manifold}, Trans. Amer. Math. Soc. \textbf{80} (1955),
  528--542. \MR{MR0084825 (18,930a)}

\bibitem{minsurf}
Kirill Krasnov and Jean-Marc Schlenker, \emph{Minimal surfaces and particles in
  3-manifolds}, Geom. Dedicata \textbf{126} (2007), 187--254. \MR{MR2328927}

\bibitem{volume}
Kirill Krasnov and Jean-Marc Schlenker, \emph{On the renormalized volume of
  hyperbolic 3-manifolds.}, Comm. Math. Phys. \textbf{279} (2008), 637--668,
  {math.DG/0607081}.

\bibitem{lecuire}
Cyril Lecuire, \emph{Plissage des vari\'et\'es hyperboliques de dimension 3},
  Invent. Math. \textbf{164} (2006), no.~1, 85--141. \MR{MR2207784
  (2006m:57029)}

\bibitem{conebend}
Cyril Lecuire and Jean-Marc Schlenker, \emph{The convex core of quasifuchsian
  manifolds with particles}, In preparation, 2008.

\bibitem{lubotzky-magid}
Alexander Lubotzky and Andy~R. Magid, \emph{Varieties of representations of
  finitely generated groups}, Mem. Amer. Math. Soc. \textbf{58} (1985),
  no.~336, xi+117. \MR{MR818915 (87c:20021)}

\bibitem{luo-tian}
Feng Luo and Gang Tian, \emph{Liouville equation and spherical convex
  polytopes}, Proc. Amer. Math. Soc. \textbf{116} (1992), no.~4, 1119--1129.

\bibitem{MM}
Y.~Matsushima and S.~Murakami, \emph{Vector bundle valued harmonic forms and
  automorphic forms on a symmetric riemannian manifold}, Annals of Math.
  \textbf{78} (1963), 365--416.

\bibitem{mess}
Geoffrey Mess, \emph{Lorentz spacetimes of constant curvature}, Geom. Dedicata
  \textbf{126} (2007), 3--45. \MR{MR2328921}

\bibitem{montcouquiol-these}
Gr\'egoire Montcouquiol, \emph{D\'eformations de m\'etriques {Einstein} sur des
  vari\'et\'es \`a singularit\'es coniques}, Ph.D. thesis, Universit\'e
  Toulouse III, 2005.

\bibitem{montcouquiol1}
\bysame, \emph{Rigidit\'e infinit\'esimale de c\^ones-vari\'et\'es {Einstein}
  \`a courbure n\'egative}, preprint math.DG/0503195, 2005.

\bibitem{montcouquiol2}
\bysame, \emph{{Deformations Einstein infinitesimales de cones-varietes
  hyperboliques}}, 2006, arXiv:math.DG/0603514.

\bibitem{horo}
Jean-Marc Schlenker, \emph{Hypersurfaces in {$H\sp n$} and the space of its
  horospheres}, Geom. Funct. Anal. \textbf{12} (2002), no.~2, 395--435.
  \MR{MR1911666 (2003d:53108)}

\bibitem{hmcb}
\bysame, \emph{Hyperbolic manifolds with convex boundary}, Invent. Math.
  \textbf{163} (2006), no.~1, 109--169. \MR{MR2208419 (2006m:57023)}

\bibitem{takhtajan-teo}
Leon~A. Takhtajan and Lee-Peng Teo, \emph{Liouville action and
  {W}eil-{P}etersson metric on deformation spaces, global {K}leinian
  reciprocity and holography}, Comm. Math. Phys. \textbf{239} (2003), no.~1-2,
  183--240. \MR{MR1997440 (2005c:32021)}

\bibitem{thurston-notes}
William~P. Thurston, \emph{Three-dimensional geometry and topology.},
  Originally notes of lectures at Princeton University, 1979. Recent version
  available on http://www.msri.org/publications/books/gt3m/, 1980.

\bibitem{weil}
A.~Weil, \emph{On discrete subgroups of {Lie} groups}, Annals of Math.
  \textbf{72} (1960), no.~1, 369--384.

\bibitem{weil:remarks}
Andr{\'e} Weil, \emph{Remarks on the cohomology of groups}, Ann. of Math. (2)
  \textbf{80} (1964), 149--157. \MR{MR0169956 (30 \#199)}

\bibitem{weiss-local}
Hartmut Weiss, \emph{Local rigidity of 3-dimensional cone-manifolds}, J.
  Differential Geom. \textbf{71} (2005), no.~3, 437--506. \MR{MR2198808
  (2007e:53048)}

\bibitem{weiss-global}
\bysame, \emph{Global rigidity of 3-dimensional cone-manifolds}, J.
  Differential Geom. \textbf{76} (2007), no.~3, 495--523. \MR{MR2331529
  (2009b:53063)}

\end{thebibliography}

\def\cprime{$'$}
\providecommand{\bysame}{\leavevmode\hbox to3em{\hrulefill}\thinspace}
\providecommand{\MR}{\relax\ifhmode\unskip\space\fi MR }
\providecommand{\MRhref}[2]{%
  \href{http://www.ams.org/mathscinet-getitem?mr=#1}{#2}
}
\providecommand{\href}[2]{#2}

\end{document}